\documentclass[a4paper,12pt,intlimits,oneside]{amsart}
\usepackage{amsmath}
\usepackage{amsthm}
\usepackage{latexsym}
\usepackage{amssymb}
\usepackage{xcolor}
\numberwithin{figure}{section}
\def\R{{\mathbb R}}
\def\C{{\mathbb C}}

\def\Z{{\mathbb Z}}
\def\S{{\mathbb S}}

\def\N{{\mathbb N}}
\def\s{\vskip 0.25cm\noindent}
\def\build#1_#2^#3{\mathrel{
\mathop{\kern 0pt#1}\limits_{#2}^{#3}}}
\def\td_#1,#2{\mathrel{\mathop{\build\longrightarrow_{#1\rightarrow #2}^{}}}}
\newtheorem{theorem}{Theorem}[section]
\newtheorem{corollary}{Corollary}
\newtheorem{proposition}{Proposition}
\newtheorem{Lemma}{Lemma}
\newtheorem{remark}{Remark}
\newtheorem{definition}{Definition}
\begin{document}
\title{The cubic Szeg\"o equation}
\author{Patrick G\'erard}
\address{Universit\'e Paris-Sud XI, Laboratoire de Math\'ematiques
d'Orsay, CNRS, UMR 8628} \email{{\tt Patrick.Gerard@math.u-psud.fr}}

\author[S. Grellier]{Sandrine Grellier}
\address{MAPMO-UMR 6628,
D\'epartement de Math\'ematiques, Universit\'e d'Orleans, 45067
Orl\'eans Cedex 2, France} \email{{\tt
Sandrine.Grellier@univ-orleans.fr}}

\subjclass[2000]{ 35B15, 37K10, 47B35}

\date{June 24, 2009}

\keywords{Nonlinear Schr\"odinger equations, Integrable Hamiltonian
systems, Lax pairs, Hankel operators} \maketitle
\renewcommand{\abstractname}{R\'esum\'e}
\begin{abstract}
On consid\`ere l'\'equation hamiltonienne suivante sur l'espace de
Hardy du cercle
$$i\partial _tu=\Pi(|u|^2u)\ ,$$
o\`u $\Pi$ d\'esigne le projecteur de Szeg\"o. Cette \'equation est
un cas mod\`ele d'\'equation sans aucune propri\' et\' e dispersive.
On \'etablit qu'elle admet une paire de Lax et une infinit\'e de
lois de conservation en involution, et qu'elle peut \^etre
approch\'ee par une suite de syst\`emes hamiltoniens de dimension
finie compl\`etement int\'egrables. N\'eanmoins, on met en
\'evidence des ph\'enom\`enes d'instabilit\'e illustrant la
d\'eg\'en\'erescence de cette structure compl\`etement int\'
egrable. Enfin, on caract\'erise les ondes progressives de ce
syst\`eme.
\end{abstract}
\renewcommand{\abstractname}{Abstract}
\begin{abstract}
We consider  the following Hamiltonian equation on the $L^2$ Hardy
space on the circle,
$$i\partial _tu=\Pi(|u|^2u)\ ,$$
where $\Pi $ is the Szeg\"o projector. This equation can be seen as
a toy model for totally non dispersive evolution equations. We
display a Lax pair structure for this equation. We prove that it
admits  an infinite sequence of conservation laws in involution, and
that it can be approximated by a sequence of finite dimensional
 completely integrable Hamiltonian systems.
We establish several instability phenomena illustrating the
degeneracy of this completely integrable structure. We also classify
the traveling waves for this system.
\end{abstract}
\thanks {The authors would like to thank S.~Alinhac, L.~Baratchart,
T.~Kappeler, S.~Kuksin, J.~Leblond, W.~Strauss and M.~Zworski for
valuable discussions, and, for their hospitality, the IMPA in Rio de
Janeiro, the Chennai Mathematical Institute and the CIRM in Luminy,
where part of this work was done. They also acknowledge the supports
of the following ANR projects :  EDP dispersives
(ANR-07-BLAN-0250-01) for the first author, and AHPI
(ANR-07-BLAN-0247-01) for the second author.}
\begin{section}{Introduction}
\subsection{Motivation}
This work can be seen as a continuation of a series of papers due to
N.~Burq, N.~Tzvetkov and the first author \cite{BGT0,
BGT1,BGT2,BGT3} --- see also \cite{PG} for a survey--- , devoted to
 the influence of the geometry of a Riemannian manifold $M$ onto the qualitative
 properties of solutions to the nonlinear
Schr\"odinger equation,
\begin{equation}\label{schrod} i\partial
_tu+\Delta u=|u|^2u\ , (t,x)\in \R \times M\ . \end{equation} The
usual strategy for finding global solutions to the Cauchy problem is
to solve locally in time in the energy space $H^1\cap L^4$  using a
fixed point argument and
 then to globalize in time, by means of conservation of energy and of $L^2$ norm.
As a corollary of the work of Burq, G\'erard, Tzvetkov --- see
\cite{BGT2}, remark 2.12 p.205,  one obtains, whatever the geometry
is, the following general result. If there exists a smooth local in
time flow map on the Sobolev space $H^s(M)$, then the following
Strichartz--type estimate must hold,
\begin{equation}\label{smoothflow}
\| {\rm e}^{it\Delta }f\| _{L^4([0,1] \times M)}\lesssim \| f\|
_{H^{s/2}(M)}\ .
\end{equation}
This inequality is valid for instance if $M=\R ^d, d=1,2,3,4$ and
$\Delta$ is the Euclidean Laplacian, where $s$ is given by the
scaling formula
$$s=\max (0, \frac d2-1)\ .$$
In \cite{BGT1, BGT2}, it is observed that, on the two-dimensional
sphere, the infimum of the numbers $s$ such that (\ref{smoothflow})
holds is $1/4$, hence is larger than the regularity given by the
latter formula. This can be interpreted as a lack of dispersion
properties for the spherical geometry. It is therefore natural to
ask whether there exist some geometries for which these dispersion
properties totally disappear. Such an example arises in
sub-Riemannian geometry, more precisely for radial solutions of the
Schr\"odinger equation associated to the  sub-Laplacian on the
Heisenberg group, as observed in \cite{GG}, where part of the
results of this paper are announced. Here we  present a more
elementary example of such a situation. Let us choose $M=\R ^2
_{x,y}$ and replace the Laplacian by the Grushin operator
$G:=\partial _x^2+x^2\partial _y^2$, so that our equation is
\begin{equation}\label{eq}
i\partial _tu+\partial _x^2u+x^2\partial _y^2u=\vert u\vert ^2u\ \ .
\end{equation}
Notice that this equation enjoys the following scaling invariance :
if $u(t,x,y)$ is a solution, then
$$\lambda u(\lambda ^2t,\lambda x,\lambda ^2y)$$
is also a solution. In this context it is natural to replace the
standard Sobolev space $H^s(M)$  by the Grushin Sobolev space $H^s_G
(M)$, defined as the domain of $\sqrt {(-G)^{s}}\ .$ Observe that
the above scaling transformation leaves invariant the homogeneous
norm of $H^{1/2}_G(M)$, which suggests that equation (\ref{eq}) is
{\it subcritical} with respect to the energy regularity $H^1_G(M)$.
However, we are going to see that (\ref{smoothflow}) cannot hold if
$s<\frac 32$, which means that no smooth flow can exist on the
energy space, hence equation (\ref{eq}) should rather be regarded as
{\it supercritical} with respect to the energy regularity. In fact,
the critical regularity $s_c=\frac 32$ is  the regularity which
corresponds to the Sobolev embedding in $M$, since $x$ has
homogeneity $1$ and $y$ has homogeneity $2$. This is an illustration
of a total lack of dispersion for equation (\ref{eq}). \s
 The
justification is as follows. Notice that $u=e^{itG}f$ can be
explicitly described by using the Fourier transform in the $y$
variable, and by making an expansion along the Hermite functions
$h_m$ in the $x$ variable, leading to the representation
$$u(t,x,y)=(2\pi )^{-1/2}\sum _{m=0}^\infty \int _\R {\rm e}^{-it(2m+1)|\eta
|+iy\eta}\hat f_m(\eta )h_m(\sqrt {|\eta |}x)\, d\eta \ ,$$ with
$$\| f\|^2 _{H^{s/2}_G}=\sum _{m=0}^\infty \int _\R (1+(2m+1)|\eta |)^{s/2}|\hat f_m(\eta
)|^2\, \frac{d\eta }{\sqrt {|\eta |}}.$$ Let us focus onto data
concentrated on modes $m=0, \eta \sim N^2$, specifically
$$f(x,y)= \frac 1{\sqrt {N}}\int _0^\infty {\rm e}^{iy\eta -\eta \frac {x^2}{2}-\frac \eta {N^2}}\, d\eta =N^{\frac 32}F(Nx,N^2y) $$
with
$$F(x,y):= \frac 1{1+\frac {x^2}2 -iy}\ .$$
Then the above formula for $u$ gives
$$u(t,x,y)=f(x,y-t)\ ,$$
so that
$$\| u\| _{L^4([0,1]\times \R ^2_{x,y})}=N^{3/4}\| F\| _{L^4}\ .$$ Since $\| f\| _{H^{s/2}_G}\simeq N^{s/2}$ as $N\rightarrow \infty $,
this proves the claim. \s Let us study the {\sl structure of the
nonlinear evolution problem (\ref{eq}}). Denote by $V^\pm _m$ the
space of functions of the form
$$v^\pm _m(x,y)=\int _0^\infty e^{\pm i\eta y}g(\eta )h_m( \sqrt \eta x)\, d\eta \ ,\
\int _0^\infty \eta ^{-1/2}\vert g(\eta )\vert ^2\, d\eta <\infty \
,$$ so that we have the orthogonal decomposition
$$L^2(M)=\oplus _{\pm }\oplus _{\ m = 0}^\infty V^\pm _m\ ,\ G_{|V^\pm _m}=\pm
i(2m+1)\partial _y\ .$$ Denote by
 $\Pi _m ^\pm : L^2(M)\rightarrow
V_m^\pm $ the orthogonal projection. Expanding the solution as
$$u=\sum _\pm \sum _{m=0}^\infty u_m^\pm \ ,\ u_m^\pm =\Pi _m^\pm u\ ,$$
 the equation reads as a system of coupled transport equations,
\begin{equation}\label{system}
i(\partial _t\pm (2m+1)\partial _y)u_m=\Pi _m^\pm (|u|^2u)\ .
\end{equation}
Therefore a better understanding of equation (\ref{eq}) requires to
study the interaction between the nonlinearity $|u|^2u$ and the
projectors $\Pi _m^\pm $.  Notice that similar interactions arise in
the literature, see for instance \cite{Nier} in the study of the
Lowest Landau Level for Bose-Einstein condensates, or
\cite{BGTstrauss} in the study of critical high frequency regimes of
NLS on the sphere. Other examples can be found in the introduction
of \cite{GG}. The present paper is devoted to a toy model for this
kind of interaction.
\subsection{A toy model : the cubic Szeg\"o equation}
Let $$\S^1= \{z\in \C , |z|=1\} $$ be the unit circle in the complex
plane. If $u$ is a distribution on $\S^1$, $u\in \mathcal D'(\S^1)$,
then $u$ admits a Fourier expansion
 in the distributional sense
$$u=\sum _{k\in \Z}\hat u(k)e^{ik\theta }\ .$$
For every subspace $E$ of ${\mathcal D}'(\S^1)$, we denote by $E_+$
the subspace
$$E_+=\{ u\in E\ ; \ \forall k<0, \hat u(k)=0\} \ .$$
In particular, $L^2_+$ is the Hardy space of $L^2$ functions which
extend to the unit disc $\{ |z|<1\} $ as holomorphic functions,
$$u(z)=\sum _{k=0}^\infty \hat u(k)z^k\ ,\ \sum _{k=0}^\infty |\hat u(k)|^2<+\infty \ .$$
 Let us endow
$L^2(\S^1)$ with the scalar product
$$(u|v):=\int _{\S^1}u\overline v\, \frac{d\theta }{2\pi }\ ,$$
and denote by $\Pi :L^2(\S^1)\rightarrow L^2_+(\S^1)$ be the
orthogonal projector  on $L^2_+(\S^1)$, the so-called Szeg\"o
projector,
$$\Pi \left (\sum _{k\in \Z}\hat u(k)e^{ik\theta }\right )=\sum _{k\ge 0}\hat u(k)e^{ik\theta }
. $$ We consider the following evolution equation on $L^2_+(\S^1)$,
\begin{equation}\label{S}
i\partial _tu=\Pi (|u|^2u)\ . \end{equation} This equation, that we
decided to call {\sl the cubic Szeg\"o equation}, is the simplest
one which displays  interaction between a cubic nonlinearity and a
Calderon-Zygmund projector. It is also an infinite dimensional
Hamiltonian system on $L^2_+(\S^1)$, as we shall now see.
\subsection{The Hamiltonian formalism}
We endow $L^2_+(\S^1)$ with the symplectic form
$$\omega (u,v)=4\,{\rm Im}(u|v)\ .$$
Given a real valued function $F$ defined on a dense subspace
${\mathcal D}$ of $L^2_+(\S^1)$, we shall say that $F$ admits a
Hamiltonian vector field if there exists a mapping
$$X_F:{\mathcal D}\rightarrow L^2_+(\S^1)$$
such that, for every $h\in {\mathcal D}$,
$$\frac{F(u+th)-F(u)}{t}\td _t,0 \omega (h, X_F(u))\ .$$
Of course,  this property is often strengthened as differentiability
of  $F$ for some norm on ${\mathcal D}$ (see Kuksin \cite{K} for a
general setting in scales of Hilbert spaces). A Hamiltonian curve
associated to $F$ is a solution $u=u(t)$ of
$$\dot u=X_F(u)\ ,$$
and, given two functions $F, G$ on ${\mathcal D}$ admitting
Hamiltonian vector fields, the Poisson bracket of $F, G$ is defined
on ${\mathcal D}$ by
$$\{ F,G\}(u)=\omega (X_F(u), X_G(u))\ .$$
For example, the function
$$E(u)=\int _{\S^1}|u|^4\, \frac{d\theta}{2\pi }\ ,$$
defined on $L^4_+(\S^1)$, admits on $H^s_+(\S^1)$, $s>\frac 12$,
the Hamiltonian vector field
$$X_E(u)=-i\Pi (|u|^2u)\ ,$$
which defines a smooth vector field on $H^s_+$, so that equation
(\ref{S}) is the equation of Hamiltonian curves for $E$.
 From
this structure, the equation $(S)$ inherits the formal conservation
law $E(u)=E(u(0))$. The invariance by translation and by
multiplication by complex numbers of modulus $1$ gives two other
formal conservation laws,
$$Q(u):=\int _{\S^1}|u|^2\, \frac{d\theta }{2\pi
}=\|u\|_{L^2}^2\quad ,\quad M(u):=(Du|u),\ D:=-i\partial _\theta
=z\partial _z\ .$$ Equivalently, these conservation laws mean that
we have the following cancellations for the Poisson brackets,
$$\{ E, Q\} =\{ E, M\} =0\ ,$$
which can be recovered in view of the explicit expressions of the
Hamiltonian vector fields,
$$X_Q(u)=-\frac i2u\ ,\ X_M(u)=-\frac i2 Du\ .$$
Finally, these expressions also imply that
$$\{ Q, M\} =0\ .$$
\subsection{Main results}
>From the previous conservation laws, we shall show --- see section
\ref{sCauchy}--- that (\ref{S}) defines a continuous flow on
$H^{1/2}_+$. The main results of this paper are based on an
unexpected property of this flow, namely that it admits a Lax pair,
as the KdV flow (see Lax \cite{L}) or the one dimensional cubic
Schr\"odinger flow (see Zakharov-Shabat \cite{Z}). More precisely,
for every $u\in H^{1/2}_+$, we define (see {\it e.g.} Peller
\cite{P}, Nikolskii \cite{N}), the Hankel operator of symbol $u$ by
$$H_u(h)=\Pi (u\overline h)\ ,\ h\in L^2_+\ .$$
It is well known that $H_u$ is a Hilbert-Schmidt operator, which is
symmetric with respect to the real part of the scalar product on
$L^2_+$. Our basic result is roughly the following --- see section
\ref{laxpair} for a more precise statement.
\begin{theorem}
There exists a mapping $u\mapsto B_u$, valued into skew--symmetric
operators on $L^2_+$, such that $u$ is a solution of (\ref{S}) if
and only if
$$\frac d{dt} H_u=[B_u,H_u]\ .$$
\end{theorem}
\noindent As a consequence, if $u$ is a solution of (\ref{S}),
$H_{u(t)}$ is unitarily equivalent to $H_{u(0)}$. From this
observation, we infer many new properties of the dynamics of
(\ref{S}), including an infinite  sequence $(J_{2n})_{n\ge 1}$ of
conservation laws in evolution.
 We also prove  the approximation of equation (\ref{S}) by finite dimensional
completely integrable Hamiltonian systems --- see sections
\ref{M(N)} and \ref{Hierarchy}.
\begin{theorem}\label{W}
For every positive integer $D$, there exists a complex submanifold
$W(D)$ of $H^{1/2}_+$ of dimension $D$, such that
\begin{enumerate}
\item $W(D)$ is invariant by the flow of (\ref{S}).
\item The flow of (\ref{S}) is a completely integrable Hamiltonian flow
on $W(D)$ in the Liouville sense.
\end{enumerate}
Moreover,  the union of the manifolds $W(D)$, $D\ge 1$, is dense in
$H^{1/2}_+$.
\end{theorem}
\noindent In Theorem \ref{W} above, complete integrability in the
Liouville sense means,  according to Arnold \cite{Ar}, that for
generic Cauchy data in $W(D)$, the evolution is quasi-periodic on a
Lagrangian torus. In fact, $W(D)$ is a manifold of rational
functions on the complex plane, with no poles in the unit disc. For
instance, $W(3)$ consists of functions $u$ given by
$$u(z)=\frac{az+b}{1-pz}$$
with $a\in \C \setminus \{ 0\} , b\in \C ,$ and $p$ in the open unit
disc. In this particular case, we solve (\ref{S}) explicitly  in
section \ref{Mtilde1}, and we deduce the following large time
behavior of $H^s$ norms of the solutions.
\begin{theorem}\label{WTW(3)}
Every solution $u$ of (\ref{S}) on $W(3)$ satisfies
$$\forall s>\frac 12\ ,\ \sup _{t\in \R}\|u(t)\| _{H^s}<+\infty .$$
However, there exists a family $(u_0^\varepsilon )_{\varepsilon >0}$
of Cauchy data in $W(3)$, which converges in $W(3)$ for the
$C^\infty (\S ^1)$ topology as $\varepsilon \rightarrow 0$, such
that the corresponding solutions $u^\varepsilon $ satisfy
$$\forall \varepsilon >0, \exists t^\varepsilon >0 : \forall s>\frac 12\ ,
\| u^\varepsilon (t^\varepsilon )\| _{H^s}\td _\varepsilon, 0
+\infty .$$
\end{theorem}
\noindent The second statement of Theorem \ref{WTW(3)} is to be
compared with the recent result by
Colliander-Keel-Staffilani-Takaoka-Tao \cite{CKSTT}, who proved a
similar behavior for cubic NLS on the two-dimensional torus. Notice
that, as shown by our result, this behavior does not imply the
existence of an unbounded trajectory in $H^s$, and that it can occur
for completely integrable systems. This phenomenon shows that the
conservation laws of equation (\ref{S}) do not control the high
energy Sobolev norms. However, let us mention  that the boundedness
of the trajectories in $H^s$ is a generic property on all the
manifolds $W(D)$, as we prove in section \ref{largetime}. The
boundedness in $H^s$ of the trajectory for all data in $H^s$ for
large $s$, is an interesting open problem. \s Finally, in section
\ref{TravelingWaves} we characterize traveling waves for (\ref{S}).
In view of the two-dimensional symmetry group associated to $Q$ and
$M$, these traveling waves are defined as follows.
 \begin{definition}\label{introtrav}
 A solution $u$ of (\ref{S})  is said to be a traveling wave
 if there exists $\omega,c\in \R$ such  that $$u(t,z)=e^{-i\omega t}u(0, e^{-ict}z)$$ for
  every $t\in \R $. We shall call $\omega $ the {\it pulsation} of $u$, and $c$ the velocity of $u$.
\end{definition}
\noindent Notice that the equation for traveling waves is the
following nonlinear equation,
$$cDu+\omega u=\Pi (|u|^2u)\ .$$
In section \ref{TravelingWaves}, using the Lax pair structure and a
precise spectral analysis of the corresponding selfadjoint
operators, we describe all the solutions of this equation.
\begin{theorem}
The initial data $u_0\in H^{1/2}_+$ of traveling waves for (\ref{S})
are given by \begin{eqnarray*} u_0(z)=\begin{cases}
\displaystyle{\alpha \prod _{j=1}^N\frac{z-\overline p_j}{1-p_jz}\
{\rm for}\
\alpha \in \C, |p_j|<1, N\ge 1,\ {\rm if}\  c=0}\ ,\\
\displaystyle{\alpha \frac{z^\ell }{1-p^Nz^N}\ ,\  {\rm for}\ \alpha
\in \C, N\geq 1, 0\leq \ell \leq N-1\ ,\ {\rm if}\ c\ne 0.}
\end{cases}\end{eqnarray*}
\end{theorem}
\noindent The question of orbital stability of these traveling
waves, in the sense of Grillakis-Shatah-Strauss \cite{GSS2}, is of
course very natural. We only have partial answers to this question,
namely in the case $N=1$ of the above theorem:
\begin{enumerate}
\item For $|p|<1$, the stationary wave corresponding to
$$u_0(z)=\frac{z-\overline p}{1-pz}$$
is orbitally unstable --- see section \ref{Mtilde1}.
\item For $|p|<1$,
the stationary wave corresponding to
$$u_0(z)=\frac{1}{1-pz}$$
is orbitally stable --- see section \ref{M(1)}. In fact, we show
that this data is a ground state of the variational equation which
characterizes the traveling waves.
\end{enumerate}
We close this introduction by mentioning two  natural open problems,
on which we hope to come back in a future work. The first one is  to
obtain a complete solution of equation (\ref{S}) by solving inverse
spectral problems for Hankel operators,
 describing explicit action
angle coordinates for (\ref{S}), as it is done in \cite{KP} for the
KdV equation.  The second one is of course to transfer at least part
of the structure found here for attacking the open problem of global
smooth solutions to
 the nonlinear Schr\"odinger equation associated to the Grushin operator, which was the starting point of this paper,
 and to other evolution problems on the same type \cite{GG}, for instance on
 the Heisenberg group.
\end{section}
\begin{section}{The Cauchy problem}\label{sCauchy}
In this section, we solve the Cauchy problem for the cubic Szeg\"o
equation, for sufficiently smooth data. We close the section by a
remark about the smoothness of the flow map. Further results
concerning uniform continuity of this flow map for weaker topologies
can be found in Section \ref{M(1)}, Proposition \ref{instabM(1)}.
\begin{theorem}\label{Cauchy}
Given $u_0\in H^{1/2}_+(\S^1)$, there exists a unique solution $u\in
C(\R ,H^{1/2}_+(\S^1))$ of (\ref{S}) such that $u(0)=u_0$. For every
$T>0$, the mapping $u_0\in H^{1/2}_+\mapsto u\in C([-T,T],
H^{1/2}_+)$ is continuous. Moreover, if $u_0\in H^s_+(\S^1)$ for
some $s>\frac 12$, then $u\in C(\R ,H^s_+(\S^1))$.
\end{theorem}
\begin{proof}
Assume first that $s>1/2$. Since the vector field $X_E$ is smooth on
$H^s_+$,  it is easy to solve $(S)$ locally in time . More
precisely, one has to solve the integral equation
\begin{equation}\label{I}
u(t)= u_0-i\int_0^t \Pi(|u|^2u) dt'.
\end{equation}
The corresponding operator is well defined on $H^s_+(\S^1)$ since
$$\|\Pi(|u|^2u)\|_{H^s}\le \||u|^2u\|_{H^s}\leq C\|u\|_{L^\infty}^2\|u\|_{H^s}\le C'\| u\| _{H^s}^3 .$$
This allows to use a fixed point argument on a  small time interval,
and yields a time interval of existence $[-T,T]$ where $T$ is
bounded from below if $\| u_0\| _{H^s}$ is bounded. \s Next we show
that the $H^s$-norm of this unique solution remains bounded on any
time interval, so that this solution is  global. To that purpose, we
make use of the conservation of $Q$ and $M$, and of the following
observation,
 \begin{equation}\label{H1/2}
M(u)+Q(u)=\sum _{k\ge 0}(k+1)|\hat u(k)|^2=\| u\| _{H^{1/2}}^2\ .
\end{equation}
So far, we have only observed that $M$ and $Q$ are formally
conserved. In fact, it  is straightforward to prove this
conservation for sufficiently smooth solutions, and finally we get
them for $H^s$ solutions, $s>1/2$, by approximation.\s  We combine
the conservation of the $H^{1/2}$ norm  with the following
Brezis-Gallou\"et type estimate (see \cite{BG}),
$$\| u\| _{L^\infty }\leq C_s\| u\| _{H^{1/2}}\left [\log \left (2+\frac{\| u\| _{H^s}}{\| u\| _{H^{1/2}}}\right )\right ]^{\frac 12}\ .$$
A proof of this  estimate is recalled in Appendix 1. We infer, for
$t\ge 0$,
\begin{eqnarray*}
\|u\|_{H^s}&\le& \|u_0\|_{H^s}+\int_0^t\|\Pi(|u|^2u)\|_{H^s} dt' \le
\|u_0\|_{H^s}+C\int_0^t \|u\|_{L^\infty}^2\|u\|_{H^s} dt'\cr &\le&
\|u_0\|_{H^s}+B\int_0^t  \| u_0\| _{H^{1/2}}^2\left [\log \left
(2+\frac{\| u\| _{H^s}}{\| u_0\| _{H^{1/2}}}\right )\right
]\|u\|_{H^s}dt'\ .
\end{eqnarray*}
If we set $f(t):=\|u\|_{H^s}/\| u_0\| _{H^{1/2}}$, we obtain
$$f(t)\le  f(0)+A\int_0^t \left [\log  (2+f(t'))\right ] f(t')dt'\ .$$
so that, by a non linear Gronwall lemma, $f$ does not blow up in
finite time , \begin{equation}\label{doublexp} 2+f(t)\leq
(2+f(0))^{{\rm e}^{At}}\ . \end{equation} This completes the proof
for $s>1/2$. \s Let us turn to the case $s=1/2$. The proof of global
existence of weak solutions is standard. Let us recall it briefly.
Given $u_0\in H^{1/2}_+$, approximate it by a sequence $(u_0^n)$ of
elements in $H^{s}_+$, $s>1/2$. Consider the sequence $(u_n)$ of
solutions of (\ref{S}) in $C(\R, H^s_+)$ corresponding to these
initial data. In view of (\ref{H1/2}), the $H^{1/2}$ norm of
$u_n(t)$ remains bounded for any $t\in\R$, and consequently
$\partial _tu_n(t)$ remains bounded in, say, $L^2$. Hence
 there exists a subsequence of $u_n(t)$
converging weakly to $u(t)$ in $H^{1/2}$, locally uniformly in $t$.
By the Rellich theorem, $u_n(t)$ converges strongly to $u(t)$ in
$L^p$ for every $p<\infty $, and  it is easy to check that such a
function $u$ is a weak solution of (\ref{S}).\s Next, let us prove
the uniqueness, which follows from an argument first introduced by
Yudovich in the case of the $2D$ Euler equation and used by
Vladimirov in \cite{V}, and Ogawa in \cite{O}.  It is based on the
fact that functions in $H^{1/2}(\S^1)$ satisfy the Trudinger-type
inequality,
\begin{equation}\label{E}
\forall p\in [1, \infty [\, , \, \| u\| _{L^p}\leq C\, \sqrt {p}\,
\| u\| _{H^{1/2}}\quad \end{equation} We postpone the proof of this
estimate to Appendix 2. Let $u$ and $\tilde u$ be two solutions of
(\ref{S}) belonging to $C_w(\R ,H^{1/2}_+)$ with $u(0)=\tilde u(0)$.
Set $g(t):=\|u(t)-\tilde u(t)\|_{L^2}^2$ so that $g$ is $C^1$ and
vanishes at the origin. Introduce a large number $p>2$ and compute
\begin{eqnarray*}
|g'(t)| &=& 2\, \left|{\rm Im} \left (\,  (u(t)-\tilde u(t))\, |\,
\Pi(|u|^2u-|\tilde u|^2\tilde u)\, \right )\right|\cr
&\le&C_1\int_{\S^1} |u-\tilde u|^2(|u|^2+|\tilde u|^2)d\theta \cr
&\le&C_1'\int_{\S^1} |u-\tilde u|^{2(1-\frac 1p)}(|u|^2+|\tilde
u|^2)^{1+\frac 1p}d\theta \cr &\le&C_2 \|u-\tilde
u\|_{L^{2}}^{2(1-\frac 1p)}(\|u\|_{L^{2(p+1)}}^{2(1+\frac
1p)}+\|\tilde u\|_{L^{2(p+1)}}^{2(1+\frac 1p)})\, \cr &\le& B\,  p\,
 g(t)^{1-\frac 1p}\ .
\end{eqnarray*}
This implies $$g(t)\le (Bt)^{p}\ .$$
 The right hand side of the latter inequality goes to zero as $p$ goes to infinity
for any $t<1/B$. This  proves the uniqueness of the Cauchy problem.
\s It remains to prove that the weak solution $u$ is strongly
continuous in time with values in $H^{1/2}$, and that it depends
continuously on the Cauchy data $u_0$. First, by weak convergence,
we have $\|u(t)\|_{H^{1/2}}\le \|u_0\|_{H^{1/2}}$ for any $t\in\R$.
By reversing time and using uniqueness, one obtains the converse
inequality for any $t\in \R$
--- solve the Cauchy problem with initial data $u(t)$. Hence the
$H^{1/2}$ norm is preserved by the flow on $H^{1/2}_+$. Since $u$ is
weakly continuous with respect to $t$ and since $H^{1/2}_+$ is a
Hilbert space, this completes the proof of the strong continuity of
$u$. The continuity of the flow map can be proved similarly.
\end{proof}

\begin{remark} For $s>1/2$,
the contraction mapping  argument used to construct the solution $u$
classically allows to prove that the flow map $u_0\mapsto u(t)$ is
Lipschitz continuous on bounded subsets of $H^s$ and that it is
smooth. \s On the opposite, the flow defined on $H^{1/2}_+(\S^1)$ is
not smooth
--- in fact it is not $C^3$ near $0$.
Here is the argument. If $\Phi _t$ is the flow map, a simple
expansion shows that, for $h\in H^s_+, s>\frac 12$,
$$d^3\Phi _t(0)(h,h,h)=-6it \Pi (|h|^2h)\ .$$
Hence the fact that $\Phi _1$ is $C^3$ on  a neighborhood of $0$ in
$H^{1/2}_+$ is in contradiction with the existence of $h\in
H^{1/2}_+$ such that
 $\Pi(|h|^2h)$ does not belong to  $H^{1/2}_+$.
 As a simple computation shows, an example of such a function $h$ is given by $h_\alpha=f^\alpha$ where $f(z)=-\frac{\log(1-z)}z$
 and $\frac 16<\alpha <\frac 12$.
\end{remark}

\end{section}
\begin{section}{A Lax pair for the cubic Szeg\"o equation.}\label{laxpair}
In this section, we show that the cubic Szeg\"o equation (\ref{S})
enjoys a very rich property, namely it admits a Lax pair in the
sense of Lax \cite{L}. As a preliminary step, we introduce relevant
operators on the Hardy space $L^2_+(\S^1)$ (see Nikolskii \cite{N}
and Peller \cite{P} for general references).\s
 Given $u\in
H^{1/2}_+(\S^1)$, the {\it Hankel operator} of symbol $u$ is defined
by
$$H_u(h)=\Pi (u\overline h)\
.$$ Notice that $H_u$ is $\C $-antilinear, and is always a symmetric
operator  with respect to the real scalar product ${\rm Re}(u|v)$.
In fact, it satisfies the identity
$$(H_u(h_1)|h_2)=(H_u(h_2)|h_1)\ .$$
 Consequently,  $H_u^2$ is $\C $-linear, selfadjoint and nonnegative.
 Moreover, $H_u$ is given in
terms of Fourier coefficients by
$$\widehat {H_u(h)}(k)=\sum _{\ell \geq 0}\hat u(k+\ell )\overline {\hat h(\ell )}\ .$$
Consequently, we have
$$\widehat{H_u^2(h)}(k)=\sum _{j\leq 0}c_{kj}h_j\ ,\ c_{kj}:=\sum _{\ell \geq 0}\hat u(k+\ell)\overline {\hat u_(j+\ell)}\ .$$
In particular, \begin{equation}\label{HS} Tr(H_u^2)=\sum _{k\geq
0}c_{kk}=\sum _{\ell \geq 0}(\ell +1)|\hat u(\ell )|^2=M(u) +Q(u)\ ,
\end{equation}
hence $H_u$ is a Hilbert-Schmidt operator. \s Given $b\in L^\infty
(\S^1)$, the  {\it Toeplitz operator} of symbol $b$ is defined by
$$T_b(h)=\Pi (bh)\ .$$
The  operator $T_b $ is of course $\C $ -linear, and is selfadjoint
for the Hermitian scalar product (hence symmetric for the real
scalar product)  as soon as $b$ is real valued.
\begin{theorem}\label{thlaxpair}
Let $u\in C(\R ,H^s(\S^1))$ for some $s>\frac 12$. The cubic Szeg\"o
equation
$$i\partial _tu=\Pi (|u|^2u)$$
is equivalent to the fact that the Hankel operator $H_u$ satisfies
the evolution equation
\begin{equation}\label{evolop}
\frac{d}{dt}H_u=[B_u,H_u]
\end{equation}
where \begin{equation}\label{B} B_u=\frac i2 H_u^2-iT_{|u|^2}
\end{equation} is a skew-symmetric operator. In other words, the pair
$(H_u,B_u)$ is a Lax pair for the cubic Szeg\"o equation.
\end{theorem}
\begin{proof}
Firstly, we establish the following identity,
\begin{equation}\label{Rio}
H_{\Pi (|u|^2u)}=T_{|u|^2}H_u+H_uT_{|u|^2}-H_u^3 .
\end{equation}
Given $h\in L^2_+$, we have
$$H_{\Pi (|u|^2u)}(h)=\Pi (\Pi (|u|^2u)\overline h)=\Pi
(|u|^2u\overline h)$$ since $\Pi ((1-\Pi )(b)\overline h)=0$ for
every $b$. Then
$$\Pi
(|u|^2u\overline h)=\Pi (|u|^2\Pi (u\overline h))+\Pi (|u|^2(1-\Pi
)(u\overline h)),$$ and we observe that
$$\Pi (|u|^2\Pi (u\overline h))=T_{|u|^2}H_u(h)\ ,$$
while
$$\Pi (|u|^2(1-\Pi
)(u\overline h))=H_u\left (u\overline{(1-\Pi)(u\overline h)}\right
)\ .$$ It remains to notice that, since
$u\overline{(1-\Pi)(u\overline h)}\in L^2_+$, \begin{eqnarray*}
u\overline{(1-\Pi)(u\overline h)}&=&\Pi \left
(u\overline{(1-\Pi)(u\overline h)}\right )\\&=&\Pi (|u|^2h)-\Pi
\left (u\overline {\Pi (u\overline h)}\right )=
T_{|u|^2}(h)-H_u^2(h)\ .
\end{eqnarray*}
This completes the proof of (\ref{Rio}). Now we just observe that
(\ref{S}) is equivalent to
$$\frac d{dt}H_u=-iH_{\Pi (|u|^2u)}=[B_u,H_u]$$
since $H_u$ is antilinear.
\end{proof}
As a consequence of Theorem \ref{thlaxpair}, the cubic Szeg\"o
equation admits an infinite number of conservation laws. Indeed,
from (\ref{evolop}), we classically observe that, denoting by $U
(t)$ the solution of the operator equation
$$\frac{d}{dt}U=B_u\, U, U(0)=I\ ,$$
the operator $U(t)$ is unitary for every $t$, and
$$U(t)^*H_{u(t)}U(t)=H_{u(0)}\ .$$
 In other words, we have the following property.
\begin{corollary}\label{isospec}
Let $u$ be a solution of (\ref{S}) with initial value $u_0\in H^s_+,
s>1/2$. The family of Hankel operators $(H_{u(t)})_{t\in \R}$ is
isospectral to $H_{u_0}$.
\end{corollary}

\medskip

Let us state some consequences of this isospectrality. First, we
recall some basic properties of Hankel operators (see \cite{N},
\cite{P} for proofs). It is well known from a theorem by Nehari
\cite{Ne}  that the operator norm of $H_u$ is equivalent to $\| u\|
_{BMO} +\| u\| _{L^2}$, which is therefore essentially conserved by
the flow. Moreover, a theorem by Peller states that, for $p<\infty
$, the Schatten norm $[{\rm Tr}(|H_u|^p)]^{1/p}$
 is equivalent to the norm of $u$ in  the Besov space $B^{1/p}_{p,p}$, which is therefore uniformly bounded
 for all time if it is finite at $t=0$.
 Notice that the particular case $p=2$ was already observed in
 (\ref{HS}), giving again the conservation of $M(u)+Q(u)$.
Another example of a conserved quantity is of course the trace norm
$Tr(|H_u|)\ ,$ which, as stated before,  is equivalent to the Besov
$B^1_{1,1}$ norm of $u$ (or to the $L^1$-norm of $u''$ with respect
to the area measure in the disc). This observation leads to a
significant improvement of the large time estimate (\ref{doublexp})
for the high Sobolev norms of the solution of (\ref{S}) derived from
the proof of Theorem \ref{Cauchy}.
\begin{corollary}\label{sobolev}
Assume $u_0\in H^{s} _+$ for some $s>1$. Then we have the following
estimates,
\begin{eqnarray*}
\sup _ {t\in \R}\| u(t)\| _{L^\infty }&\leq & C\| u_0\| _{H^s}\ ,\cr
\| u(t)\| _{H^s}&\leq & C \| u_0\| _{H^s}\, {\rm e}^{C\| u_0\|
_{H^s}|t|}\  .
\end{eqnarray*}
\end{corollary}
\begin{proof} Since $H^s\subset B^{1}_{1,1}$ as soon as $s>1$, the trace norm
of $H_{u_0}$ is finite, hence the $B^{1}_{1,1}$ norm of $u(t)$ is
uniformly bounded. Since $B^{1}_{1,1}\subset L^\infty $, this proves
the first assertion. The second one is then a simple consequence of
the standard Gronwall lemma. \end{proof} We will return to the large
time behavior of solutions of (\ref{S}) in sections \ref{Mtilde1}
and \ref{largetime}. At this stage, it is natural to find a way to
recover  other known conservation laws, namely $Q$ and $E$. In fact,
we are going to find them as two particular cases of an  infinite
sequence of conservation laws, which will play an important role in
the sequel.
\begin{corollary}\label{Jn}
For every $u\in H^{1/2}_+$, for every positive integer $n$, set
$$J_n(u)=(H_u^n(1)|1)\ .$$
If $u\in C(\R ,H^{1/2}_+)$ solves (\ref{S}), we have, for every
positive integer $k$,
$$\frac{d}{dt}J_{2k}(u)=0\ ,\ i\frac{d}{dt}J_{2k-1}(u)=J_{2k+1}(u)\
.$$
\end{corollary}
\begin{proof}
We may assume that $u_0\in H^s$ for $s>1/2$, since the general case
follows by density and the continuity properties of the flow map on
$H^{1/2}_+$. Coming back to (\ref{evolop}), we observe that
$$B_u(1)=\frac i2 H_u^2(1)-iT_{|u|^2}(1)=-\frac i2H_u^2(1)\ .$$
Consequently, since $H_u^{2k}$ is $\C $-linear and $B_u$ is skew
symmetric,
\begin{eqnarray*}
\frac d{dt}(H_u^{2k}(1)|1)&=&( [B_u,H_u^{2k}](1),|1)\\
&=&-(H_u^{2k}(1)|B_u(1))-(H_u^{2k}B_u(1)|1)\\
&=& -\frac i2(H_u^{2k+2}(1)|1)+\frac i2(H_u^{2k+2}(1)|1) =0\ .
\end{eqnarray*}
The second identity is obtained similarly, observing that
$H_u^{2k-1}$ is $\C$ -antilinear,
\begin{eqnarray*}
i\frac d{dt}(H_u^{2k-1}(1)|1)&=&i( [B_u,H_u^{2k-1}](1),|1)\\
&=&-i(H_u^{2k-1}(1)|B_u(1))-i(H_u^{2k-1}B_u(1)|1)\\
&=& \frac 12(H_u^{2k+1}(1)|1)+\frac 12(H_u^{2k+1}(1)|1)
=J_{2k+1}(u)\ .
\end{eqnarray*}
\end{proof}
The conservation of $Q$ and $E$ is recovered by observing that
\begin{eqnarray*}
J_2(u)&=&(H_u^{2}(1)|1)_{L^2}=\| u\| _{L^2}^2=Q(u)\ ,\\
J_4(u)&=&(H_u^4(1)|1)_{L^2}=\| H_u^2(1)\| _{L^2}^2=\|\Pi (|u|^2)\|
_{L^2}^2=\frac {E(u)+Q(u)^2}{2}\ .
\end{eqnarray*}
In section \ref{Hierarchy}, we will prove that the conservation laws
$J_{2k}$ are in involution, and that that their differentials
satisfy some generic independence.
\end{section}

\begin{section}{Invariant finite dimensional submanifolds}\label{M(N)}
In this section, we introduce finite dimensional submanifolds of
$L^2_+$ which are invariant by the flow of the cubic Szeg\"o
equation. Elements of these manifolds turn out to be rational
functions of the variable $z$, with no poles in the unit disc. In
what follows, $\C _D[z]$ denotes the class of complex polynomials of
degree at most $D$, and $d(A)$ denotes the degree of a polynomial
$A$.
\subsection{The manifold ${\mathcal M}(N)$}
\begin{definition} Let $N$ be a positive integer.
We denote by ${\mathcal M}(N)$ the set of rational functions $u$ of
the form
$$u(z)=\frac {A(z)}{B(z)}\ ,$$
with $A\in \C_{N-1} [z] $, $B\in \C_N [z],$ $B(0)=1,$ $d(A)=N-1$ or
$d(B)= N$, $A$ and $B$ have no common factors,  and $B(z)\ne 0$ if
$|z|\leq 1$.
\end{definition}
Notice that ${\mathcal M}(N)$ is included in  $H^s_+$ for every $s$.
It is elementary to check that ${\mathcal M}(N)$ is a
$2N$-dimensional complex submanifold of $L^2_+$, and that its
tangent space at $u=A/B$ is
$$T_u{\mathcal M}(N)= \frac{\C_{2N-1}[z]}{B^2}\ .$$
 A theorem by Kronecker states
 that ${\mathcal M}(N)$ is exactly the set of symbols $u$ such
 that $H_u$ is of  rank $N$.
 For the convenience of the reader, we give an elementary proof of
 this result in Appendix 3.
 In view of Corollary \ref{isospec}, we infer the following result, which can also be checked directly, using some elementary linear algebra.
 \begin{theorem}
 Let $u_0\in{\mathcal M}(N)$ and $u$ be the solution of (\ref{S}) with $u(0)=u_0$.
 Then, for every $t\in\R$, $u(t)$ belongs to ${\mathcal M}(N)$.
 In other words, the submanifolds ${\mathcal M}(N)$ are
 invariant under the flow of the cubic Szeg\"o equation.
 \end{theorem}
\noindent In the notation of Theorem \ref{W} of the introduction,
the manifold ${\mathcal M}(N)$ is $W(2N)$. Since ${\mathcal M}(N)$
is finite dimensional, equation (\ref{S}) on ${\mathcal M}(N)$ is
reduced to a system of ordinary differential equations, which we now
describe in the main coordinate patch of ${\mathcal M}(N)$. A
generic point in ${\mathcal M}(N)$ is given by
$$u=\sum_{j=1}^N\frac{\alpha_j}{1-p_jz},$$
where the $p_j$'s are  pairwise distinct and belong to the unit
disc. Then, in the coordinates $(\alpha _j, p_j)_{1\leq j\leq N}$,
(\ref{S}) reads
\begin{eqnarray}\displaystyle\label{eqN}
\begin{aligned}\begin{cases}
i\dot \alpha _j&=\sum _k\frac{\alpha _j^2\overline {\alpha _k}}{(1-p_j\overline {p_k})^2}+2\sum _k\sum _{\ell \ne j}\frac{\alpha _j\overline {\alpha _k}\alpha _\ell p_j}{(p_j-p_\ell)(1-p_j\overline{p_k})}\ , \\
i\dot p_j&=\sum _k\frac {\alpha _j\overline {\alpha
_k}}{1-p_j\overline{p_k}}\; p_j\ ,
\end{cases}\end{aligned}
\end{eqnarray}
In particular, the conservation laws $Q$, $M$, $E$  read
\begin{eqnarray*}
Q&=&\sum _{j,k}\frac{\alpha _j\overline \alpha _k}{1-p_j\overline
p_k}\ ,\ M=\sum _{j,k}\frac{\alpha _jp_j\overline \alpha _k\overline
p_k}{(1-p_j\overline p_k)^2}\ ,\hskip
4cm  \\
E&=&\sum _{j,k,l,m}\frac{\alpha _j\overline \alpha _k\alpha
_l\overline \alpha _m(1-p_j\overline p_kp_l\overline
p_m)}{(1-p_j\overline p_k)(1-p_j\overline p_m)(1-p_l\overline
p_k)(1-p_l\overline p_m)}\ .
\end{eqnarray*}
In view of the second part of system (\ref{eqN}), we notice an
additional conservation law, \begin{equation}\label{clS}
S=|p_1\cdots p_N|^2\ .
\end{equation}
In the next subsection, we give an intrinsic interpretation of $S$
and we establish further properties which will be useful in the
sequel.
\subsection{The Blaschke product associated to $u\in {\mathcal M}(N)$}\label{blaschke}
Given $u\in H^{1/2}_+$, it is elementary to check from
$$H_u(h)=\Pi (u\overline h)$$
 that $\ker H_u$ is a
closed subspace of $H_u$ invariant by the shift $h\mapsto zh$.
According to the Beurling Theorem \cite{R}, there exists $\varphi
\in L^2_+$,  such that $|\varphi |^2=1$ on $\S^1$ and
$$\ker H_u=\varphi L^2_+\ .$$
Let us characterize such a generator $\varphi $ if $u=A/B \in
{\mathcal M}(N)$. Set
$$B(z)=\prod _{j=1}^N (1-p_jz)\ ,$$
where the $p_j$'s are complex numbers in the open unit disc, with
possible repetitions. We define the {\it Blaschke product}
associated to $u$ by
$$b(z)=\prod _{j=1}^N\frac{z-\overline p_j}{1-p_jz}\ .$$
and we claim that \begin{equation}\label{beurling} \ker H_u=bL^2_+\
.
\end{equation}
Indeed, if $u=A/B$ as in definition \ref{M(N)}, the equation $\Pi
(u\overline h)=0$ means exactly that there exists $g\in L^2_+$ such
that
$$z^{N-1}\overline A\left (\frac 1z\right )h(z)=g(z)\prod
_{j=1}^N(z-\overline p_j)\ .$$ On the other hand, the assumptions on
$A$, $B$ imply that the polynomials $z^{N-1}\overline A\left (\frac
1z\right )$ and $\prod _{j=1}^N(z-\overline p_j)$ have no common
factor. Consequently, $\ker H_u$ consists of those $h\in L^2_+$
which are divisible by $\prod _{j=1}^N(z-\overline p_j)$, which is
equivalent to $h\in bL^2_+$. \s Let us make the connection with the
distinguished vector $1$. Since ${\rm Im}(H_u)$ is finite
dimensional and since $H_u$ is symmetric, we have
$${\rm Im}(H_u)=(\ker (H_u))^\perp \ .$$
In particular, ${\rm Im}(H_u)$ is a space of rational functions,
whose general description is provided in Appendix 3. We denote by
$P_u$ the orthogonal projector on ${\rm Im}(H_u)$.
\begin{proposition}\label{v}
We have
$$1-P_u(1)=(-1)^Np_1\cdots p_N\, b\ .$$
\end{proposition}
\begin{proof}
Set $v=1-P_u(1)$. From Appendix 3, $1\in {\rm Im}(H_u)$ if and only
if one of the $p_j$'s is $0$. Since the claimed identity is trivial
in this case, we may assume that $p_j\ne 0$ for every $j$. Then,
from Appendix 3, all the functions in ${\rm Im}(H_u)$ tend to $0$ at
infinity, hence $v(z)$ tends to $1$ at infinity. Since
$$v(z)=h(z) b(z)\ ,$$
where $h$ is a polynomial, we conclude that
$$h(z)=(-1)^Np_1\cdots p_N\ .$$
\end{proof}
As a consequence, we obtain the following interpretation of the
conservation law $S$ introduced in the previous subsection,
$$S:=|p_1\cdots p_N|^2={\rm dist}(u,\ker H_u)^2\ .$$
Indeed, ${\rm dist}(u,\ker H_u)^2=\| 1-P_u(1)\| _{L^2} ^2$ and
$|b|^2=1$ on $\S^1$, hence we even have $|1-P_u(1)|^2=S$ on $\S^1$.
In fact, we can derive a more general evolution law for the whole
quantity $v=1-P_u(1)$.
\begin{proposition}\label{evolv}
Let $u$ be a solution of (\ref{S}) on ${\mathcal M}(N)$. Then
$v:=1-P_u(1)$ satisfies, on $\S^1$,
$$i\partial _tv=|u|^2v\ .$$
\end{proposition}
\begin{proof}
Notice that $v$ is the orthogonal projection of $1$ onto $\ker H_u$,
which reads, in terms of the functional calculus of the selfadjoint
operator $H_u^2$,
$$v={\bf 1}_{\{ 0\}} (H_u^2)(1)\ .$$
Consequently, by Theorem \ref{thlaxpair},
$$\partial _tv=[B_u,{\bf 1}_{\{ 0\}} (H_u^2)](1)\ .$$
Since
$$B_u=-iT_{|u|^2}+\frac i2H_u^2\ ,\ B_u(1)=-\frac i2H_u^2(1)\ ,$$
we get
$$i\partial _tv=T_{|u|^2}v\ .$$
The following lemma implies that $T_{|u|^2}v=|u|^2v$ and therefore
completes the proof.
\begin{Lemma}\label{ubarh}
If $u\in L^\infty _+\cap H^{1/2}_+$ and $h\in \ker H_u$, then
$\overline uh\in zL^2_+$.
\end{Lemma}
Indeed, for every $k\ge 0$, we have , in $L^2(\S^1)$,
$$(\overline uh|\overline z^k)=(z^k|u\overline h)=(z^k|\Pi
(u\overline h))=(z^k|H_u(h))=0\ .$$
\end{proof}
As a consequence of the above proposition, let us deduce an
evolution law for the Blaschke product $b$ if $S\ne 0$. In this
case, $v$ does not vanish on the circle, and we can write, at each
point of $\S^1$,
$$|u|^2=i\frac{\partial _tv}{v}=i\frac{\partial _t(p_1\cdots
p_N)}{p_1\cdots p_N}+i\frac{\partial _tb}{b}\ .$$ Let us take the
average of both sides on $\S^1$. Since
$$\frac{\partial _tb}{b}=\sum _{j=1}^N\left (-\frac{\partial _t\overline
p_j}{z-\overline p_j}+\frac{z\partial _tp_j}{1-p_jz}\right )\ ,$$ a
direct calculation yields
$$\int _{\S^1}\frac{\partial _tb}{b}\, \frac{dz}{2i\pi z}=0\ ,$$
and therefore
\begin{equation}\label{evolprod} Q=i\frac{\partial
_t(p_1\cdots p_N)}{p_1\cdots p_N}\ .
\end{equation}
Coming back to Proposition \ref{evolv}, we infer
\begin{equation}\label{evolb}
i\partial _tb=(|u|^2-Q)b\ .
\end{equation}
Equation (\ref{evolb}) in fact holds without assuming $S\ne 0$. This
can be shown by approximation in ${\mathcal M}(N)$. However, we
shall give a different proof in the next subsection, which is
devoted to the flow on the subset $\{ S=0\} $ of ${\mathcal M}(N)$.
\subsection{The manifold $\tilde {\mathcal M}(N-1)$}\label{tilde}
 Denote by $\tilde{ \mathcal M} (N-1)$ the subset of $\mathcal M(N)$
defined by the equation $S=0$.  The rational functions in
$\tilde{\mathcal M}(N-1)$ are the elements of ${\mathcal M}(N)$ with
a numerator of degree exactly equal to $N-1$ and a denominator of
degree at most $N-1$, therefore $\tilde {\mathcal M}(N-1)$ is a
complex hypersurface of ${\mathcal M}(N)$, and its tangent space at
$u=A/B$ is
$$T_u\tilde {\mathcal M}(N-1)=\frac{\C _{2N-2}[z]}{B^2}\ .$$
 As $S$ is invariant under the flow,
we get that $\tilde{ \mathcal M}(N-1)$ is invariant under the flow.
In the notation of Theorem \ref{W} of the introduction, $\tilde
{\mathcal M}(N-1)$ is $W(2N-1)$. On this submanifold, generic points
are described as
$$u=\sum _{j=1}^{N-1}\frac{\alpha _j}{1-p_jz}+\alpha_{N}\ ,$$
where the $p_j$'s are as before in the open unit disc, pairwise
distincts and different from $0$. The generic evolution is system
(\ref{eqN}) with $p_N=0$. From this explicit system, we notice that
the trivial conservation law $S$ is replaced  by  $$\tilde S=\left
|\alpha _{N}\prod _{j=1}^{N-1}p_j\right |^2\ .$$ As in the previous
section, we shall now give a more intrinsic interpretation of the
new conservation law $\tilde S$. \s Since $1\in {\rm Im}(H_u)={\rm
Im}(H_u^2)$, there exists a unique $w\in {\rm Im}(H_u)$ such that
$$H_u(w)=1\ .$$
Write
$$u=\frac AB\ ,\ B(z)=\prod _{j=1}^{N-1}(1-p_jz)\ ,\
A(z)=az^{N-1}+\sum _{j<N-1}a_jz^j\ ,$$ with $a\ne 0$. The associated
Blaschke product now reads
$$b(z)=z\prod _{j=1}^{N-1}\frac{z-\overline p_j}{1-p_jz}:=z\tilde
b(z)\ .$$ Notice that, from the description of ${\rm Im}(H_u)$
provided in Appendix 3, $\tilde b\in {\rm Im}(H_u)$. From the
elementary identity
$$H_u(zh)=\overline z(H_u(h)-(u|h))\ ,$$
we infer
$$H_u(\tilde b)=(u|\tilde b)\ .$$
Then an explicit calculation gives
$$(\tilde b|u)=\int _{\S^1}\frac{z^{N-1}\overline A(1/z)}{\prod
_j(1-p_jz)}\frac{dz}{2i\pi z}=\overline a\ .$$ Therefore we have
proved
\begin{proposition}\label{w}
$$w(z)=\frac{\tilde b(z)}{\overline a}=\frac{b(z)}{\overline az}\ .$$
\end{proposition}
We conclude this subsection by deriving an evolution law for $w$.
\begin{proposition}\label{evolw}
Let $u$ be a solution of (\ref{S}) on $\tilde {\mathcal M}(N-1)$.
Then the preimage $w$ of $1$ in ${\rm Im}(H_u)$ satisfies, on
$\S^1$,
$$i\partial _tw=|u|^2w\ .$$
\end{proposition}
\begin{proof}
The proof is very similar to the one of Proposition \ref{evolv}.
Firstly, we express $w$ by means of  the functional calculus of the
selfadjoint operator $H_u^2$,
$$w=f(H_u^2)H_u(1)\ ,\ f(\lambda ):=\frac{{\bf 1}_{]0,\infty [}(\lambda)}{\lambda }\ .$$
Consequently, by Theorem \ref{thlaxpair},
$$\partial _tw=[B_u,f(H_u^2)H_u](1)\ .$$
Since
$$B_u=-iT_{|u|^2}+\frac i2H_u^2\ ,\ B_u(1)=-\frac i2H_u^2(1)\ ,$$
we get
$$i\partial _tw=T_{|u|^2}w\ .$$
The following lemma implies that $T_{|u|^2}w=|u|^2w$ and therefore
completes the proof.
\begin{Lemma}\label{ubarw}
If $u\in \tilde {\mathcal M}(N-1)$,  then $\overline uw\in L^2_+$.
\end{Lemma}
Indeed, for every $k\ge 1$, we have , in $L^2(\S^1)$,
$$(\overline uw|\overline z^k)=(z^k|u\overline w)=(z^k|\Pi
(u\overline w))=(z^k|H_u(w))=(z^k|1)=0\ .$$
\end{proof}
As a consequence of Proposition \ref{evolw}, we infer that $\| w\|
_{L^2}^2$ is a conservation law. In the case of a generic element of
${\mathcal M}(N-1)$,
$$u=\sum _{j=1}^{N-1}\frac{\alpha _j}{1-p_jz}+\alpha _N\ ,$$
we have
$$a=(-1)^{N-1}p_1\cdots p_{N-1}\alpha _N\ ,$$
thus we get the interpretation of $\tilde S$ as
$$\tilde S=|a|^2=\frac{1}{\| w\| _{L^2}^2}\ .$$
Finally, as in the previous subsection, Proposition \ref{evolw}
leads to an evolution law for the coefficient $a$ itself and for
$b$. Indeed, taking the average on $\S^1$ of
$$|u|^2=i\frac{\partial _tw}{w}=-i\frac{\partial _t(\overline a)
}{\overline a}+i\frac{\partial _t\tilde b}{\tilde b} =
-i\frac{\partial _t(\overline a) }{\overline a}+i\frac{\partial
_tb}{b}\ ,
$$
 we obtain
 \begin{equation}\label{evola}
 i\partial _ta=Qa\ ,
 \end{equation}
and, coming back to the equation on $w$, we eventually deduce the
evolution of $b$ (\ref{evolb}), in the whole generality. \s In the
next two sections, we study the particular cases of $\mathcal M(1)$
and of $\tilde{\mathcal M}(1)$ in more detail.
\end{section}

\begin{section}{The case of $\mathcal M(1)$}\label{M(1)}
Elements of ${\mathcal M}(1)$ are
\begin{equation}\label{phiap}
\varphi _{\alpha ,p}(z)=\frac{\alpha }{1-pz}\ ,\ \alpha \ne 0\ ,\
|p|<1\ .
\end{equation}
In this particular case, the system (\ref{eqN}) reads
$$i\dot \alpha =\frac{|\alpha |^2}{(1-|p|^2)^2}\alpha \ ,\ i\dot p=\frac{|\alpha |^2}{1-|p|^2}p\ ,$$
which is solved as
$$\alpha (t)=\alpha (0)\, {\rm e}^{-i\omega t}\ ,\ p(t)=p(0)\, {\rm e}^{-ict}\ ,\ \omega =\frac{|\alpha (0)|^2}{(1-|p(0)|^2)^2}\
 \ c=\frac{|\alpha (0)|^2}{1-|p(0)|^2}\ .$$
Equivalently, the solution $u$ of (\ref{S}) with $u(0)=\varphi
_{\alpha ,p}$ is given by
$$u(t,z)={\rm e}^{-i\omega t}\, \varphi _{\alpha ,p}({\rm
e}^{-ict}z),$$ which means that $u$ is a traveling wave according to
Definition \ref{introtrav}. In section \ref{TravelingWaves}, we will
classify all such solutions. Notice that, apart from the trivial
case of constants --- $p=0$ ---, the trajectory lies in the
two-dimensional torus $\{ |\alpha |=cst\ ,\ |p|=cst \} $. We are
going to prove that this two-dimensional torus can also be seen as
the solution of a variational problem in $H^{1/2}_+$. We first state
the following lemma which is an easy consequence of the
Cauchy-Schwarz inequality.
\begin{Lemma}\label{OpPositif}
Let $A$ be a positive operator on a separable Hilbert space
${\mathcal H}$ and $e$ be an  element of $\mathcal H$ so that
$Ae\neq 0$. Then, the following inequality holds
$$\Vert Ae\Vert^2 \le ( Ae|e) \text{Tr}(A).$$
Furthermore, equality holds if and only if $A$ is of rank one.
\end{Lemma}
Applying this lemma to $A=H^2_u$ on ${\mathcal H}=L^2_+$ with $e=1$,
and using the formulae for $J_2$ and $J_4$ derived in section
\ref{laxpair}, we get the following characterization of the elements
of ${\mathcal M}(1)$, which can be seen as an analogue of M.
Weinstein's sharp Gagliardo-Nirenberg inequality \cite{W1}.
\begin{proposition}\label{CaractM(1)} For every $u\in H^{1/2}_+$,
$$E(u)\leq Q(u)(Q(u)+2M(u)),$$
i.e.
$$\| u\| _{L^4}^4\leq \| u\| _{L^2}^2(\| u\| _{L^2}^2+\| u\| _{H^{1/2}}^2)\  ,$$
with equality if and only if $u\in\mathcal M(1)$.
\end{proposition}
Let us mention that a more direct proof of Proposition
\ref{CaractM(1)} can be found in \cite{GG}. As a consequence of
Proposition \ref{CaractM(1)}, we obtain the following large time
stability of ${\mathcal M}(1)$ in $H^{1/2}_+$.
\begin{corollary}\label{stabM(1)}
Let  $a>0 \, ,\, 0<r<1,$ and
$$T(a,r)=\{ \varphi _{\alpha ,p} : |\alpha |=a\ ,\ |p|=r\}\ .$$
For every $\varepsilon >0$, there exists $\delta >0$ such that, if
$u_0\in H^{1/2}_+ $ satisfies
$$\inf _{\varphi \in T(a,r)}\| u_0-\varphi\| _{H^{1/2}}\leq \delta $$
then the solution $u$ of (\ref{S}) with $u(0)=u_0$ satisfies
$$\sup _{t\in \R}\ \inf _ {\varphi \in T(a,r)}\|u(t)-\varphi\| _{H^{1/2}}\leq \varepsilon \ .$$
\end{corollary}
\begin{proof}
By Proposition \ref{CaractM(1)} and a simple calculation of
$Q(\varphi _{\alpha ,p})$ and $E(\varphi _{\alpha ,p})$,
 $T(a,r)$ is the set of minimizers of the problem
$$\inf \{ M(u)\ ,\ u\in H^{1/2}_+\ ,\ Q(u)=q(a,r)\ ,\  E(u)=e(a,r)\} =m(a,r)\ ,$$
where
$$q(a,r):=\frac{a^2}{1-r^2}\ ,\ e(a,r):=\frac{a^4(1+r^2)}{(1-r^2)^3}\ .$$
Let $u_0^n$ be a sequence of $H^{1/2}_+$ such that
$$\inf _{\varphi \in T(a,r)}\| u_0^n-\varphi \|  _{H^{1/2}}\rightarrow 0\ .$$
Then
$$Q(u_0 ^n)\rightarrow q(a,r)\ ,\ E(u_0^n)\rightarrow e(a,r)\ ,\ M(u_0^n)\rightarrow m(a,r)$$
and by the conservation laws,
$$Q(u^n(t))\rightarrow q(a,r)\ ,\ E(u^n(t))\rightarrow e(a,r)\ ,\ M(u^n(t))\rightarrow m(a,r)$$
 uniformly in $t$.
Given any sequence $(t_n)$ of real numbers, the  sequence
$(u^n(t_n))$ is bounded in $H^{1/2}_+$, hence has  a subsequence
which converges  weakly to some $u$ in $H^{1/2}_+$, and we get, by
the weak continuity of $Q, E$ and the weak semi-continuity of $M$,
$$Q(u)=q(a,r)\ ,\ E(u)=e(a,r)\ ,\ M(u)\leq m(a,r),$$
hence finally $M(u)=m(a,r)$, which implies from Proposition
\ref{CaractM(1)} that $u\in T(a,r)$ and that $u^n(t_n)$ converges
strongly to $u$. The proof is complete.
\end{proof}
The explicit evolution  of (\ref{S}) on ${\mathcal M}(1)$ also
allows to prove the following high frequency instability result
 in $H^s_+$ for every $s<1/2$. This result means that, given a time
 $t\ne 0$, the flow map at time $t$ does not extend as  a uniformly
 continuous map on bounded subsets of $H^s_+$, $s<\frac 12$, or
 $L^4_+$ (see Tzvetkov \cite{T} for a general discussion).
\begin{proposition} \label{instabM(1)}
Let $s<\frac 12$. There exist $u_0^\varepsilon $, $\tilde
u_0^\varepsilon $ bounded  sequences in $H^s_+$ such that
$$\| u_0^\varepsilon -\tilde u_0^\varepsilon \| _{H^s}\rightarrow 0\ {\rm but}\ \forall t\ne 0,\ \liminf _{\varepsilon \rightarrow 0}\| u^\varepsilon (t)-\tilde u^\varepsilon (t)\| _{H^s}>0\ .$$
The same holds for $H^s_+$ replaced by $L^4_+$.
\end{proposition}
\begin{proof} The principle of the proof follows Birnir-Kenig-Ponce-Svansted-Vega \cite{BSV}.
 As $|p|\rightarrow 1$, one has
$$\displaystyle {\left \| \varphi_{\alpha,p}\right \|  ^2_{H^s}
=\left \| \frac{\alpha }{1-pz}\right \|  ^2_{H^s}\sim \frac{|\alpha
|^2}{(1-|p|^2)^{1+2s}}\ .}$$ Choose
$$u_0^\varepsilon=\varphi_{\varepsilon^{s+\frac 12},\sqrt {1-\varepsilon}}\ ,\
\tilde u_0^\varepsilon=\varphi_{\varepsilon ^{s+\frac 12}(1+\delta )
,\sqrt{1-\varepsilon}}\ , $$ with $\delta \rightarrow 0$ so that $\|
u_0^\varepsilon -\tilde u_0^\varepsilon \| _{H^s}\rightarrow 0$. By
the previous computations, we get $u(t,{\rm e}^{i\theta})={\rm
e}^{-i\omega t} u_0({\rm e}^{i(\theta-ct)})$ and $\tilde u(t,{\rm
e}^{i\theta})= {\rm e}^{-i\tilde \omega t} \tilde u_0({\rm
e}^{i(\theta-\tilde ct)})$ where
$\tilde c-c=\varepsilon^{2s}\delta(2+\delta)$. Choose
$\varepsilon\rightarrow 0$ so that $\delta \varepsilon
^{2s-1}\rightarrow \infty $. It implies in particular that
$\displaystyle{\frac {\tilde c-c}{\varepsilon }\rightarrow \infty
.}$

We claim that, for any $t>0$,
$$ \| u^\varepsilon (t)-\tilde u^\varepsilon (t)\| _{H^s}^2=\| u^\varepsilon (t)\| _{H^s}^2+\|\tilde u^\varepsilon (t)\| _{H^s}^2+o(1)$$
 as $\varepsilon$ goes to zero. In other words, the scalar product in $H^s$ of $u^\varepsilon (t)$ and $\tilde u^\varepsilon (t)$ is $o(1)$.
 The result will follow since $\| u^\varepsilon (t)\| _{H^s}\simeq \| \tilde u^\varepsilon (t)\| _{H^s}\simeq 1$.

We have
\begin{eqnarray*}
|\langle u^\varepsilon (t), \tilde u^\varepsilon
(t)\rangle_{H^s}|&=&|\sum_k(1+|k|^2)^{s} \widehat{
u^\varepsilon}(t,k)\cdot\overline
{\widehat{\tilde u^\varepsilon}(t,k)}|\\
&=&|\sum_k (1+|k|^2)^{s}e^{-ik(c-\tilde c)t} \widehat u_0^\varepsilon(k)\cdot\overline{\widehat{\tilde u_0^\varepsilon}(k)}|\\
&=&\varepsilon^{2s+1}(1+\delta)|\sum_k (1+|k|^2)^{s}e^{-ik(c-\tilde c)t} (1-\varepsilon)^k|\\
&\simeq& \frac{\varepsilon^{2s+1}}{|1-(1-\varepsilon)e^{-i(c-\tilde c)t}|^{1+2s}}\\
&\simeq&\left(\frac{\varepsilon}{|c-\tilde
c|t}\right)^{1+2s}=o(1)t^{-(1+2s)}\ .
\end{eqnarray*}
The proof for $L^4_+$ is similar, observing that
$$ \| \varphi_{\alpha,p}\|_{L^4}^4=\frac{|\alpha
|^4(1+|p|^2)}{(1-|p|^2)^3}\ .$$ Choose the same functions
$u_0^\varepsilon $ and $\tilde u_0^\varepsilon$ as above, with
$s=\frac 14$, and $\delta $ going to $0$ such that $\delta
\varepsilon ^{-1/2} \rightarrow \infty $. In view of the explicit
expression, $$|u_0^\varepsilon (e^{i\theta })|^4=\frac{\varepsilon
^3}{(2-\varepsilon-2\sqrt {1-\varepsilon } \cos \theta )^2}\ ,$$ one
easily checks that, if $R^\varepsilon \rightarrow \infty $,
$$\int _{R^\varepsilon \varepsilon <|\theta -ct|<\pi }|u^\varepsilon (t,e^{i\theta
})|^4\, d\theta \rightarrow 0 \ ,\ \int _{R^\varepsilon \varepsilon
<|\theta -\tilde ct|<\pi }|\tilde u^\varepsilon (t,e^{i\theta
})|^4\, d\theta \rightarrow 0 \  .$$ Let us choose $R^\varepsilon $
such that
$$R^\varepsilon <\! < \frac{\tilde c-c}{\varepsilon }\ .$$
Then we claim that, for $t\ne 0$,
$$ \| u^\varepsilon (t)-\tilde u^\varepsilon (t)\| _{L^4}^4=\| u^\varepsilon (t)\| _{L^4}^4+\|\tilde u^\varepsilon (t)\| _{L^4}^4+o(1)\ .$$
Indeed, if $a+b=4$ and $a,b\ge 1$, by H\"older inequality, we have
$$\int _{{\mathbb S}^1}|u^\varepsilon (t)|^a|\tilde u^\varepsilon
(t)|^b\, d\theta =\int _{E^\varepsilon }|u^\varepsilon (t)|^a|\tilde
u^\varepsilon (t)|^b\, d\theta +o(1)\ ,$$ where
$$E^\varepsilon =\{\theta \in (-\pi, \pi )\ :\ |\theta -ct|< R^\varepsilon
\varepsilon\ ,\  |\theta -\tilde ct|< R^\varepsilon \varepsilon \} \
.$$ In view of the assumption on $R^\varepsilon $, this set is empty
for $\varepsilon $ small enough. This completes the proof.
\end{proof}
\end{section}
\begin{section}{The case of $\tilde{\mathcal M}(1)$}\label{Mtilde1}
The manifold  $\tilde{\mathcal M}(1)$ is a three-dimensional
K\"ahler manifold, on which (\ref{S}) admits  three conservation
laws in involution, which are $Q, M, E$. As we will see later, these
conservation laws are generically independent on $\tilde{\mathcal
M}(1)$, therefore the equation $(S)$ is  completely integrable on
this manifold.
 We are going to solve this system explicitly, by introducing coordinates which are close to the
  action angle coordinates provided by the Liouville theorem (see
  Arnold \cite{Ar}). Then we will establish some instability phenomena for large time.
 \begin{subsection}{The evolution on $\tilde{\mathcal M}(1)$}
Let us make some preliminary calculations. Since the rank of $H_u^2$
is $2$ if  $u\in \tilde {\mathcal M}(1)$, the Cayley-Hamilton
theorem reads
\begin{equation}\label{CH}
H_u^4-\sigma _1H_u^2+\sigma _2 P_u=0\ . \end{equation} Here, $\sigma
_1$ is the trace of $H_u^2$ so it equals $Q+M$. Let us compute
$\sigma _2$. Applying the above formula to the preimage $w\in {\rm
Im}(H_u)$ of $1$ introduced in subsection \ref{tilde}, we get
\begin{equation}\label{CHw}
H_u^3(1)-(M+Q)u+\sigma _2w=0\ .
\end{equation}
Taking the scalar product of (\ref{CHw}) with $w$, and using that
$(u|w)=(H_u(1)|w)=(H_u(w)|1)=1$, we infer
$$\sigma _2=M\tilde S\ .$$
We now apply (\ref{CH}) to $1\in {\rm Im}(H_u)$, and take the scalar
product with $1$. This yields
$$J_4=(M+Q)Q-M\tilde S$$
or
$$E=Q^2+2M(Q-\tilde S)\ .$$
Consequently, we can use $M,Q,\tilde S$ rather than $M,Q,E$ as our
three conservation laws. For future reference, we introduce the
solutions $r_\pm $ of the characteristic equation,
$$r^2-\sigma _1r+\sigma _2=0\ ,$$
given by
$$r_\pm =\frac12 \left (Q+M\pm ((Q+M)^2-4M\tilde S)^{1/2}\right )\ ,$$
and we set
$$\Omega =r_+-r_-=((M+Q)^2-4M\tilde
S)^{1/2}\ .$$
\begin{proposition}\label{evolMtilde}
Let $u_0\in \tilde {\mathcal M}(1)$, and let $u={\displaystyle
\frac{az+b}{1-pz}}$ be the corresponding solution of (\ref{S}). One
of the following two cases occurs :
\begin{itemize}
\item Either $Q=\tilde S$, and
\begin{equation}\label{statMtilde} u_0(z)=a_0\frac{z-\overline p}{1-pz}\ ,\ u(t,z)={\rm
e}^{-iQt}u_0(z)\ .
\end{equation}
\item Or $Q>\tilde S$, and the evolution of $a,b,p$ is given by
\begin{equation}\label{af}
i\dot a=Qa\ ,\ i\dot f_\pm =r_\pm f_\pm \ ,
\end{equation}
with
$$f_\pm :=r_\pm b+Ma\overline p\ .$$ In particular, $|p|^2$ satisfies
$$|p|^2=A+B\cos (\Omega t+\varphi )$$
for some constants $A, B, \varphi $,  and $|p|$ oscillates between
the following values,
\begin{equation}\label{rho}
\rho _{{\rm max}}= \frac{M^{1/2}+\tilde S^{1/2}}{(M+Q+2\sqrt
{M\tilde S})^{1/2} } \ ,\
 \rho _{{\rm min}}= \frac{|M^{1/2}-\tilde S^{1/2}|}{(M+Q-2\sqrt {M\tilde S})^{1/2}}\ .
\end{equation}
\end{itemize}
\end{proposition}
\begin{remark}
In the case (\ref{statMtilde}), the solution $u$ is called a
stationary wave. We will classify such solutions in section
\ref{TravelingWaves}.
\end{remark}
\begin{proof}
We already know that
$$i\dot a =Qa\ .$$
By corollary \ref{Jn}, we also know that $J_1=b$ and $J_3$ satisfy
$$i\dot J_1=J_3\ , i\dot J_3=J_5\ ,$$
and $J_5=\sigma _1J_3-\sigma _2J_1$ in view of  (\ref{CH}). Finally,
$J_3$ is easily obtained by taking the scalar product of (\ref{CHw})
with $1$ and using Proposition \ref{w},
$$J_3=(M+Q)J_1-M\tilde S(w|1)=(M+Q)b+Ma\overline p\ .$$
Setting
 $$f_\pm :=J_3-r_\mp J_1=r_\pm b+Ma\overline p\ ,$$
 we finally obtain the system of linear ODE  (\ref{af}).  Let us
first investigate the particular case $r_+=r_-$, which  is
equivalent to
$$(Q+M)^2-4M\tilde S=0\ .$$
Since $Q\geq \tilde S$ by the Cauchy-Schwarz inequality applied to
$u$ and $w$, we conclude that $r_+=r_-$ is equivalent to
$$M=Q=\tilde S\ .$$
Using the  Cauchy-Schwarz equality case, it is easy to check that
$Q=\tilde S$ is equivalent to the collinearity of $u_0$ and $w_0$,
namely
$$u_0(z)=a_0\frac{z-\overline p}{1-pz}\ .$$
 Finally, from this expression of $u_0$, a simple computation gives $M=Q$, hence $r_+=r_-$, and, since
 $|u_0|^2=Q$ on $\S^1$, we get
$$u(t)=u_0{\rm e}^{-iQt}\ .$$
In the case $r_+\ne r_-$, we can recover $(a(t), b(t), p(t))$ from
the variables $(a(t), f_\pm (t))$ and the conservation laws
$(M,Q,\tilde S)$. In particular,
\begin{equation}\label{Map}
Ma\overline p=\frac{r_+f_- \, -\, r_-f_+}{r_+-r_-}\ .
\end{equation}
Taking the modulus of both sides of (\ref{Map}), we conclude, in
view of the differential equations satisfied by $f_\pm $, that
$$|p|^2=A+B\cos (\Omega t+\varphi )$$
for some constants $A, B, \varphi $. Consequently, in view of
(\ref{Map}), $|p|$ oscillates between the following values,
$$\rho _{{\rm max}}=\frac{r_+|f_-|+r_-|f_+|}{M\Omega \tilde S^{1/2}}\ ,\ \rho _{{\rm min}}=\frac{\vert \, r_+|f_-|-r_-|f_+|\, \vert }
{M\Omega \tilde S^{1/2}}\ .$$ Let us compute $|f_\pm |$ in terms of
$M,Q,\tilde S$. Denote by $(e_+,e_-)$ an orthonormal basis of ${\rm
Im}H_u$ such that $\C e_\pm =\ker (H_u^2-r_\pm )$. Up to multiplying
$e_\pm $ by a suitable complex number of modulus $1$, we may assume,
using the $\C $-antilinearity of $H_u$, that
$$H_ue_\pm =\sqrt{r_\pm }e_\pm \ .$$
Then
$$1=\overline {\zeta }_+e_++\overline {\zeta }_-e_-\ ,\ \zeta _\pm :=(e_\pm |1)\ ,\ u=H_u(1)=\sqrt{r}_+\zeta _+e_++\sqrt{r_-}\zeta _-e_-\ .$$
From
\begin{eqnarray*}
1&=&|\zeta _+|^2+|\zeta _-|^2\ ,\ Q=r_+|\zeta _+|^2+r_-|\zeta _-|^2\ ,\\
J_1&=&\sqrt {r_+}\zeta _+^2+\sqrt {r_-}\zeta _-^2\ ,\ J_3=r_+\sqrt
{r_+}\zeta _+^2+r_-\sqrt {r_-}\zeta _-^2\ ,
\end{eqnarray*}
we obtain
$$|f_+|=\sqrt{r_+}(Q-r_-)\ ,\ |f_-|=\sqrt{r_-}(r_+-Q)\ ,$$
and finally (\ref{rho}), by a straightforward but tedious
calculation.
\end{proof}
In the next subsections, we shall take advantage of the oscillations
of $|p|$ in establishing instability results.
\end{subsection}
\begin{subsection}{Large time estimates of $H^s$ norms}\label{Hslarge}
Our first instability result concerns large time behavior of $H^s$
norms along trajectories of the cubic Szeg\"o equation on $\tilde
{\mathcal M}(1)$.
\begin{corollary}\label{HsMtilde}For every $u_0\in \tilde{\mathcal M}(1)$, the solution $u$ of (\ref{S}) with $u(0)=u_0$ satisfies,
for every $s>1/2$,
\begin{equation}\label{Hsfini}
\limsup _{t\rightarrow \infty }\| u(t)\|_{H^s}<+\infty \ .
\end{equation}
However, there exists a family $(u_0^\varepsilon )_{\varepsilon >0}$
of Cauchy data in $\tilde M(1)$, which converges in $\tilde M(1)$
for the $C^\infty (\S ^1)$ topology as $\varepsilon \rightarrow 0$,
and $K>0$ such that the corresponding solutions $u^\varepsilon $
satisfy \begin{equation}\label{turb} \forall \varepsilon >0, \exists
t^\varepsilon
>0 : t^\varepsilon \rightarrow \infty\ ,\ \forall s>\frac 12\ , \|
u^\varepsilon (t^\varepsilon )\| _{H^s}\ge K(t^\varepsilon )^{2s-1}\
.
\end{equation}
\end{corollary}
\begin{proof}
Writing as before
$$u(t)=\frac{a(t)z+b(t)}{1-p(t)z}\ ,$$
we already know that $a(t)$ and $b(t)$ are bounded because of the
conservation of $Q(u(t))$, so the blow up of the $H^s$ norm for
large $|t|$ would only come from the fact that $|p(t)|$ approaches
$1$. But this cannot happen since, by formula (\ref{rho}),
 $$\max _t|p(t)|=\rho _{{\rm max}} <1$$
 if $Q>\tilde S$. The other case $Q=\tilde S$ corresponds to  (\ref{statMtilde}), for which $p(t)=p(0)$.
\s Let us turn to the second assertion. Consider the family of
Cauchy data $\{u_0^\varepsilon\}_{0<\varepsilon<1}$ given by
$$u_0^\varepsilon (z)=z+\varepsilon $$
and let us look at the regime $\varepsilon \rightarrow 0$. Then a
simple computation from the previous formulae  shows that
$$|p(t)|^2=\frac{2}{4+\varepsilon^2}(1-\cos(\varepsilon
t\sqrt{4+\varepsilon^2})).$$ On the other hand,  using Fourier
expansion, we have, as $|p|$ approaches $1$,  $$\|u\|_{H^s}^2\simeq
\frac {|a+bp|^2}{(1-|p|^2)^{2s+1}}= M\frac1{(1-|p|^2)^{2s-1}}$$
since
$$M(u)= \frac{|bp+a|^2}{(1-|p|^2)^2}\ .$$
 In our particular case, $M(u)=1$ and we get, for $\displaystyle t^\varepsilon=\frac\pi{\varepsilon\sqrt{4+\varepsilon^2}}$,
 $$\|u(t^\varepsilon)\|_{H^s}^2\simeq \frac 1{(1-|p(t^\varepsilon)|^2)^{2s-1}}\simeq C (t^\varepsilon)^{2(2s-1)}.$$
This completes the proof.
\end{proof}
\begin{remark} Property (\ref{turb}) can be seen as a quantitative version of an instability property proved in \cite{CKSTT}
 for NLS on the two dimensional torus:
bounded data in $C^\infty $ may yield large solutions in $H^s$ for
large time. However, as shown by (\ref{Hsfini}), this may happen
even if the $H^s$ norms stay bounded on each individual trajectory,
and moreover in the case of a completely integrable system. Notice
that this phenomenon can occur with arbitrarily small data, since
multiplying the Cauchy data by a parameter $\delta $ amounts to
replace the solution $u(t)$ of (\ref{S})  by  $\delta u(\delta
^2t)$.
\end{remark}

\end{subsection}
\begin{subsection}{ Orbital instability of stationary waves.}
Our next instability result concerns the stationary waves in $\tilde
{\mathcal M}(1)$.
\begin{corollary}
 For each stationary wave $u_0$
of $\tilde{\mathcal M}(1)$, there exists a sequence $u_0^\varepsilon
$ which converges to $u_0$ in $C^\infty $ such that, for every $r\in
(0,1)$, there exists $t^\varepsilon $ such that the limit points in
$H^{1/2}_+$
 of $u^\varepsilon (t^\varepsilon )$ are of the form
$$v=\alpha \, \frac{z-\overline q}{1-qz}\ ,\  |\alpha |=\| u_0\| _{L^2}\ ,\  |q|=r.$$
\end{corollary}
\begin{proof}
First recall  that if $\displaystyle v=\frac{az+b}{1-pz}$ then the
conservation laws are given by $\displaystyle M=
\frac{|bp+a|^2}{(1-|p|^2)^2}$, $\displaystyle Q=
\frac{|bp+a|^2}{(1-|p|^2)} +|b|^2$ and $\tilde S=|a|^2$.\s Let
$\displaystyle u_0=a\,\frac{z-\overline p}{1-pz}$ be a stationary
wave of $\tilde{\mathcal M}(1)$. Define, for $0<\varepsilon <1$,  $$
u_0^\varepsilon=a\, \frac{(1-\varepsilon)z-\overline
p(1-\varepsilon/2)}{1-p(1+\epsilon/2)z}\ .$$ It is clear that such a
sequence converges to $u_0$ in $C^\infty $. By Proposition
\ref{evolMtilde}, for any fixed $\varepsilon$, the corresponding
solution $u^\varepsilon$ may be written as $\displaystyle
\frac{a^\varepsilon z+b^\varepsilon }{1-p^\varepsilon z}$ where
$|p^\varepsilon|$ oscillates between $\rho_{{\rm min}}^\varepsilon$
and $\rho_{{\rm max}}^\varepsilon$ given by (\ref{rho}).  Computing
these two bounds in terms of $\varepsilon$, it is easy to show that
$\rho_{{\rm max}}^\varepsilon$ tends to $1$ and $\rho_{{\rm
min}}^\varepsilon$ tends to $0$ as $\varepsilon$ goes to $0$.
Precisely, we have,
\begin{eqnarray*}
M
&=&|a|^2\frac{|(1-\varepsilon)-|p|^2(1-\varepsilon^2/4)|^2}{(1-|p|^2(1+\varepsilon/2)^2)^2}\\
&=&|a|^2\left(1-2\varepsilon\right)+\mathcal O(\varepsilon^2)\\
Q&=&|a|^2\left( \frac{|(1-\varepsilon)-|p|^2(1-\varepsilon^2/4)|^2}{(1-|p|^2(1+\varepsilon/2)^2)}+|p|^2(1-\varepsilon/2)^2\right)\\
&=&|a|^2(1-2\varepsilon)+\mathcal O(\varepsilon^2)\\
\tilde S&=&|a|^2(1-\varepsilon)^2=|a|^2(1-2\varepsilon)+\mathcal
O(\varepsilon^2).
\end{eqnarray*}
>From these estimates, we get $\Omega =\mathcal O(\varepsilon)$,
$\sqrt M-\sqrt {\tilde S}={\mathcal O}(\varepsilon ^2)$ and
\begin{eqnarray*}
\rho_{{\rm max}}^\varepsilon =1+\mathcal O(\varepsilon ^2)\ ,\
\rho_{{\rm min}}^\varepsilon= \mathcal O(\varepsilon)\ .
\end{eqnarray*}
In particular, for every $r\in (0,1)$, one can choose $t^\varepsilon
$ such that $|p^\varepsilon (t^\varepsilon )|=r.$ As the
$H^{1/2}$-norms of $u^\varepsilon(t )$ are bounded,
$u^\varepsilon(t^\varepsilon )$ has limit points in the weak
$H^{1/2}$-topology. Let $v_\infty$ be such a limit point. Since
$p^\varepsilon (t^\varepsilon )$ stays on the circle of radius $r$,
the convergence is strong and
 $v_\infty$ belongs to $\tilde{\mathcal M}(1)$. Moreover,   $Q(v_\infty)=\tilde S(v_\infty)$, hence $v_\infty $ is given by
 (\ref{statMtilde}).  This
completes the proof.

\end{proof}

We will pursue our study of large time behavior in  section
\ref{largetime}.
\end{subsection}
\end{section}
\section{Large time behavior on ${\mathcal M}(N)$}\label{largetime}
By Corollary \ref{HsMtilde}, every solution    on $\tilde{\mathcal
M}(1)$ satisfies  $$\sup _{t\in \R}\| u(t)\| _{H^s}<+\infty $$ for
$s\geq 0$. We prove that it is a generic situation on $\mathcal
M(N)$. A similar statement holds on $\tilde {\mathcal M}(N-1)$.
\begin{theorem}\label{LargeTimeB}
For every integer $N$, define
$$V_N=\{u_0\in\mathcal M(N);\; \det(J_{2(m+n)}(u_0))_{1\le m,n\le
N}=0\}\ .$$ Then $V_N$ is  a proper real analytic subvariety  of
$\mathcal M(N)$ and, for every $u_0\in \mathcal M(N)\setminus V_N$,
for every $s\ge 0$,
\begin{equation}\label{Hsb}
\sup _{t\in \R}\| u(t)\| _{H^s}<+\infty . \end{equation} In
particular, (\ref{Hsb}) holds for every $u_0$ outside a closed
subset of measure $0$. A similar statement holds on $\tilde
{\mathcal M}(N-1)$, with $$\tilde V_{N-1}:=V_N\cap \tilde {\mathcal
M}(N-1).$$
\end{theorem}
\begin{proof}
For every $u\in H^{1/2}_+$, we consider the polynomial expression
$$F_N(u)=\det(J_{2(m+n)}(u))_{1\le m,n\le N}\ .$$
 Notice that $F_N(u)=0$ if and only if the
vectors $H_u^{2k}(1)\ ,\ k=1,\cdots ,N$ are linearly dependent. In
particular, $F_N$ is identically $0$ on ${\mathcal M}(J)$ for $J<N$.
On the other hand, we shall see that $F_N$ is not identically $0$ on
${\mathcal M}(N)$. In fact, one can prove the following slightly
stronger result, which we state as a lemma for further references.
\begin{Lemma}\label{FN}
The vectors $H_u^{2k}(1), k=1,\cdots ,N,$ are generically
independent on $\tilde {\mathcal M}(N-1)$ and on ${\mathcal M}(N)$.
\end{Lemma}
\begin{proof}Indeed, if
$$u(z)=z^{N-1}+z^{N-2}\ ,$$
$u\in \tilde {\mathcal M}(N-1)$ and a simple computation shows that
the matrix of the system $1,H_u^{2}(1), H_u^4(1),\cdots ,
H_u^{2(N-1)}(1),$ in the basis $(z^j)_{0\leq j\leq N-1}$ is
triangular, hence these vectors are independent. Applying $H_u^2$,
which is one to one on ${\rm Im}(H_u)$,  the vectors $H_u^{2k}(1),
k=1,\cdots ,N,$ are independent as well, and $F_N(u)\ne 0$. Since
$\tilde {\mathcal M}(N-1)$ and ${\mathcal M}(N)$ are  connected,
this completes the proof.
\end{proof}
Theorem \ref{LargeTimeB} is then a consequence of the following
Lemma.
\begin{Lemma}\label{levelset}
If $u_0\in \mathcal M(N)\setminus V_N$, the level set
$$L_N(u_0):=\{ u\in {\mathcal M}(N): J_{2n}(u)=J_{2n}(u_0), 1\le n\le 2N\} $$
is a compact subset of ${\mathcal M}(N)$. \par \noindent If $u_0\in
\tilde {\mathcal M}(N-1)\setminus V_N$, the level set
$$\tilde L_{N-1}(u_0):=\{ u\in \tilde {\mathcal M}(N-1): J_{2n}(u)=J_{2n}(u_0), 1\le n\le 2N-1\} $$
is a compact subset of $\tilde {\mathcal M}(N-1)$.
\end{Lemma}
\begin{proof}
We just prove the statement for ${\mathcal M}(N)$. Let $u_0\in
\mathcal M(N)\setminus V_N$ and $u\in L_N(u_0)$. Let us first prove
that $M(u)=M(u_0)$. By the Cayley-Hamilton theorem applied to
$H_u^2$ on ${\rm Im}(H_u)$,
$$H_u^{2N}=\sum _{j=1}^{N}(-1)^{j-1}\sigma _j(u)H_u^{2(N-j)}\ .$$
Applying this identity to $H_u^{2p}(1)$ for $p=1,\cdots ,N$ and
taking the scalar product with $1$, we obtain  a system of $N$
linear equations in the $\sigma _j$'s,
$$J_{2(N+p)}(u)=\sum _{j=1}^{N}(-1)^{j-1}\sigma _j(u)J_{2(N+p-j)}\ ,\
1\leq p\leq N.$$ The determinant of this system is
$\det(J_{2(m+n)}(u))_{0\le m\le N-1,1\le n\le N}$, which, by the
above identity, is $(-1)^{N-1}F_N(u)/\sigma _N(u)$, hence is not
zero --- notice that $\sigma _N(u)\ne 0$, since $H_u^2$ is one to
one on ${\rm Im}(H_u)$. Solving this system, we conclude that each
$\sigma _j(u)$ is a universal function of $(J_{2n}(u))_{1\le n\le
2N}$. Since $\sigma _1=M+J_2$, this proves the claim. We infer that
every sequence of $L_N(u_0)$ is bounded in $H^{1/2}$, hence has
limit points for the weak topology of $H^{1/2}$. Let $v$ be such a
limit point. As a limit point of a sequence of ${\mathcal M}(N)$,
$v$ belongs to $\cup_{J\le N}\mathcal M(J)$. On the other hand,
since each $J_{2n}$ is continuous for  the weak topology of
$H^{1/2}$, $J_{2n}(v)=J_{2n}(u_0)$ for $n=1,\cdots ,2N$. In
particular, $F_N(v)=F_N(u_0)\ne 0$, whence $v\in {\mathcal M}(N)$
and finally $v\in L_N(u_0)$.
\end{proof}
The proof of Theorem \ref{LargeTimeB} is completed by observing that
the flow of (\ref{S}) conserves the level sets $L_N$, and that the
zeroes of the denominator of elements of a compact subset of
${\mathcal M}(N)$ do not approach the unit circle.
\end{proof}
\begin{corollary}\label{N=2}
For every $u_0\in {\mathcal M}(2)$, $s\ge 0$, (\ref{Hsb}) holds.
\end{corollary}
\begin{proof}
In view of Theorem \ref{LargeTimeB}, it is enough to consider the
case $F_2(u)=0$, which is equivalent to the collinearity of
$H_u^2(1)$ and of $P_u(1)$,
$$H_u^2(1)=\frac{Q}{1-S}P_u(1)\ .$$
If $P_u(1)=1\in {\rm Im}(H_u)$, then $|u|^2=Q$ and $u$ is a
stationary wave by Proposition \ref{stat}. If $P_u(1)\ne 1$, by
Proposition \ref{evolv}, the function  $v=1-P_u(1)$ satisfies
$$i\partial _tv=|u|^2v=
\frac{Q(1+S)}{1-S}v-\frac{Q}{1-S}v^2-\frac{QS}{1-S}\ .$$ Notice that
$S=(v|1)$ is a particular solution of this Riccati equation. Hence
we can solve it explicitly and observe that $v$ is a periodic
function of $t$ with period $2\pi /Q$. Since, by Proposition
\ref{v},
$$v(t,z)=p_1(t)p_2(t)\frac{(z-\overline
{p_1}(t))(z-\overline {p_2}(t))}{(1- p_1(t)z)(1-p_2(t)z)}\ ,$$ we
conclude that $p_1,p_2$ are periodic as well, hence cannot approach
the unit circle.
\end{proof}
\section{The Szeg\"o hierarchy}\label{Hierarchy}
 In this section, we show that the
  conservation laws $J_{2n}$ satisfy the Poisson commutation relations
 $$  \{ J_{2n}, J_{2p}\} =0\ ,$$
 and that $J_{2n}$ defines a global Hamiltonian flow for every $n$. In fact, we prove that, for every $n$, there exists a skew symmetric
 operator $B_{u,n}$ such that the pair $(H_u,B_{u,n})$ is a Lax pair for this Hamiltonian flow. The last part of the section is devoted to proving that functions $(J_{2n})_{1\le n\le 2N}$ are generically independent on ${\mathcal M}(N)$, and that functions $(J_{2n})_{1\le n\le 2N+1}$ are generically independent on $\tilde {\mathcal M}(N)$. This will complete the proof of Theorem \ref{W} in the introduction.
\begin{theorem}\label{J2n}
Let $s>\frac 12$. The map $u\mapsto J_{2n}(u)$ is smooth on $H^s_+$
and its Hamiltonian vector field is given by
\begin{equation}\label{XJ2n}
X_{J_{2n}}(u)=\frac{1}{2i}\sum
_{j=0}^{n-1}H_u^{2j}(1)H_u^{2n-2j-1}(1)\ .
\end{equation}
 Moreover,
$$
 H_{iX_{J_{2n}}(u)}=H_uA_{u,n}+A_{u,n}H_u$$
 where $A_{u,n}$ is the self adjoint operator
$$A_{u,n}(h)=\frac 14\left (\sum_{j=0}^{2n-2}H_u^{j}(1)\Pi (\overline{H_u^{2n-2-j}(1)}h)-\sum _{k=1}^{n-1}(h|H_u^{2k-1}(1))H_u^{2n-2k-1}(1)\right )\ .$$
\end{theorem}
\begin{proof}
Introduce, for $x$ real and $|x|$ small enough, the generating
functions,
$$w(x)=(1-xH_u^2)^{-1}(1)=\sum _{n=0}^\infty x^n H_u^{2n}(1)
$$
and
$$J(x,u)=(w(x)|1)=\sum _{n=0}^\infty x^n J_{2n}(u)\ .$$
We have
\begin{eqnarray*}
d_uJ(x,u).h&=&((1-xH_u^2)^{-1}x(H_uH_h+H_hH_u)(1-xH_u^2)^{-1}(1)|1)\\
&=&x[(H_uH_hw(x)|w(x))+(H_hH_uw(x)|w(x))]=2x{\rm Re}(h|w(x)H_uw(x))\\
&=&\omega (h|X(x))
\end{eqnarray*}
with
$$X(x)=\frac{x}{2i}w(x)H_uw(x)\ .$$
Identifying the coefficients of $x^n$, we get formula (\ref{XJ2n}).
The second part of the proof relies on the following lemma.
\begin{Lemma}\label{luminy}
We have the following identity,
$$H_{aH_u(a)}(h)=H_u(a)H_a(h)+H_u(a\Pi (\overline ah)-(h|a)a)\ .$$
\end{Lemma}
\begin{proof}
$$H_{aH_u(a)}(h)=\Pi(aH_u(a)\overline h)=H_u(a)H_a(h)+\Pi (H_u(a)(1-\Pi)(a\overline h))\ .$$
On the other hand,
$$(1-\Pi)(a\overline h)=\overline{\Pi (\overline ah)}-(a|h)\ .$$
The lemma follows by plugging the latter formula into the former
one.
\end{proof}
Let us complete the proof. Using the identity
$$w(x)=1+xH_u^2w(x),$$
and Lemma \ref{luminy} with $a=H_u(w)$, we get
\begin{eqnarray*}
H_{wH_u(w)}(h)&=&H_{H_u(w)}(h)+xH_{H_u(w)H_u^2(w)}(h)\\
&=&H_{H_u(w)}(h)+xH_u^2(w)H_{H_u(w)}(h)+\\
&+& xH_u\left (H_u(w)\Pi(\overline {H_u(w)}h)-(h|H_u(w))H_u(w)\right )\\
&=&wH_{H_u(w)}(h)+ xH_u\left (H_u(w)\Pi(\overline {H_u(w)}h)-(h|H_u(w))H_u(w)\right )\\
&=&w\Pi(\overline wH_uh)+xH_u\left (H_u(w)\Pi(\overline
{H_u(w)}h)-(h|H_u(w))H_u(w)\right )\ .
\end{eqnarray*}
We therefore have obtained
$$H_{wH_u(w)}=G_uH_u+H_uD_u$$
where $G_u$ and $D_u$ are the following self adjoint operators,
$$G_u(h)=w\Pi(\overline wh)\ ,\ D_u(h)=x\left (H_u(w)\Pi(\overline {H_u(w)}h)-(h|H_u(w))H_u(w)\right )\ .$$
Consequently, since $H_{wH_u(w)}$ is self   adjoint,
$$H_{wH_u(w)}=C_uH_u+H_uC_u$$
with
$$C_u=\frac 12 (G_u+D_u)\ .$$
Identifying the coefficients of $x^n$ in
$$H_{iX(x)}=\frac x2H_{w(x)H_uw(x)}\ ,$$
we infer the desired formula for $A_{u,n}$.

\end{proof}
\begin{corollary}\label{flotJ2n}
Let $s> 1$. For every $u_0\in H^s_+$, there exists a unique solution
$u\in C(\R, H^s_+)$ of the Cauchy problem
\begin{equation}
\partial _tu=X_{J_{2n}}(u)\ ,\ u(0)=u_0\ .
\end{equation}
Moreover, $u$ satisfies
\begin{equation}\label{laxpairJ2n}
\partial _tH_u=[B_{u,n},H_u]\ ,
\end{equation}
with
$$B_{u,n}(h)=\frac  {-i}4\left (\sum_{j=0}^{2n-2}H_u^{j}(1)\Pi (\overline{H_u^{2n-2-j}(1)}h)-\sum _{k=1}^{n-1}(h|H_u^{2k-1}(1))H_u^{2n-2k-1}(1)\right )\ .$$
Finally, we have the commutation identity
\begin{equation}\label{poisson}
\{ J_{2n},J_{2p}\} =0\ .
\end{equation}
\end{corollary}
\begin{proof}
The local-in-time solvability of the Cauchy problem is an easy
consequence of the fact that $H^s$ is an algebra. Moreover, to prove
global existence, it is enough to establish that the $L^\infty $
norm of $u$ does not blow up in finite time. In view of Theorem
\ref{J2n}, $u$ satisfies equation (\ref{laxpairJ2n}) on its interval
of existence. Since $B_{u,n}$ is skew symmetric, this implies that
$Tr(|H_u|)$ is conserved, and consequently, by Peller's theorem
\cite{P}, that the norm of $u$ in $B^1_{1,1}$ is bounded, and so is
the $L^\infty $ norm, whence the global existence, by an elementary
Gronwall argument. \s It remains to prove the commutation identity
(\ref{poisson}). This is equivalent to the fact that $J_{2p}$ is a
conservation law of the Hamiltonian flow of $J_{2n}$. The latter
fact is a consequence, as in section \ref{laxpair}, of  equation
(\ref{laxpairJ2n}), and of the formula
$$B_{u,n}(1)=\frac {-i}4\sum _{\ell =0}^{n-1}J_{2n-2\ell -2}H_u^{2\ell
}(1).$$
\end{proof}
We conclude this section with a complete integrability result.
\begin{corollary}\label{completint}
Let $N\geq 1$. The following properties hold.
\begin{enumerate}
\item The functions $J_{2k}\, , k=1,\cdots ,2N$ are independent in the complement of a closed subset of measure $0$ of ${\mathcal M}(N)$.
\item The functions $J_{2k}\, , k=1,\cdots ,2N+1$ are independent in the complement of a closed subset of measure $0$ of $\tilde {\mathcal M}(N)$.
\end{enumerate}
Consequently, for generic Cauchy data in ${\mathcal M}(N)$ and in
$\tilde {\mathcal M}(N)$, the solution of (\ref{S}) is
quasiperiodic.
\end{corollary}
\begin{proof}
First notice that $X_{J_{2n}}$ is tangent to ${\mathcal M}(N)$ and
to $\tilde {\mathcal M}(N)$. This can be seen either from the
explicit expression (\ref{XJ2n}) of  $X_{J_{2n}}$ compared to the
explicit description of the tangent spaces of ${\mathcal M}(N)$ and
$\tilde {\mathcal M}(N)$ in section \ref{M(N)}, or  as a consequence
of the Kronecker theorem,
$$
{\mathcal M}(N)=\{ u: rk (H_u)=N\} \ ,\  \tilde {\mathcal M}(N)=\{ u
\in {\mathcal M}(N+1): 1\in {\rm Im}(H_u)\ \} \ ,$$ compared with
the Lax pair property for the flow of $X_{J_{2n}}$ proved in
Corollary \ref{flotJ2n}. Consequently, the functions $J_{2k}$
restricted to the symplectic manifolds ${\mathcal M}(N)$ and to
$\tilde {\mathcal M}(N)$ are in involution. Therefore the second
statement of the corollary is reduced to properties $(1)$ and $(2)$.
Notice that property $(1)$ holds for $N=1$. Indeed, the linear
dependence of $J_2$ and $J_4$ at $u$ is equivalent to the fact that
$u$ is a stationary wave, which, on ${\mathcal M}(1)$, means that
$u$ is a constant. We shall prove that, for all $N$, property $(1)$
implies property $(2)$ and that property $(2)$ implies property
$(1)$ for $N+1$. This will complete the proof by induction. \s We
first prove that property $(2)$ for $N$ implies property $(1)$ for
$N+1$. We represent  the current generic point $u\in {\mathcal
M}(N+1)$ as
$$u(z)=\frac{A(z)}{B(z)}\ ,\ A\in \C _N[z],\ d(A)=N\ ,$$
with $B(z)=bz^{N+1}+\tilde B(z)$, $\tilde B\in \C_N[z]$. In this
representation, $\tilde {\mathcal M}(N)$ is characterized by the
cancellation of the holomorphic coordinate $b$. Notice that
$S=|b|^2$. Fix $u_0\in \tilde{\mathcal M}(N)$ such that the
differential form
 $$\alpha :=\bigwedge _{k=1}^{2N+1}dJ_{2k}$$
satisfies $\alpha (u_0)\ne 0$ on $T_{u_0}\tilde {\mathcal M}(N)=\ker
db(u_0)$. In  a small neighborhood $U$ of $u_0$ in ${\mathcal
M}(N+1)$,
    define 2N+1 vector fields  $Y_k, k=1,\cdots ,2N+1,$ such that, for
  every $u\in U$, $(Y_k(u))_{1\leq k\leq 2N+1}$ is a basis of $\ker (db(u))$. Since $$\alpha (u_0)(Y_1(u_0),\cdots ,Y_{2N+1}(u_0))\ne
  0,$$
  this is still true
  near $u_0$. On the other hand, since $S=|b|^2$, $dS.Y_j=0$ by construction. Hence
  \begin{eqnarray*}
  (dS\wedge \alpha )\left (b\frac{\partial }{\partial b}, Y_1,\cdots ,Y_{2N+1}\right )&=&dS\left (b\frac{\partial }{\partial b}\right )\, \alpha (Y_1,\cdots ,Y_{2N+1})\\
  &=& 2S\, \alpha (Y_1,\cdots ,Y_{2N+1})\ ,
\end{eqnarray*}
which does not cancel on $U\setminus \tilde {\mathcal M}(N)$. This
shows that the functions $S, J_{2k}, k=1,\cdots, 2N+1$ are
generically independent on ${\mathcal M}(N+1)$. In view of Lemma
\ref{FN}, we also know that the $N+1$ vectors $H_u^{2k}(1),
k=1,\cdots ,N+1$ are generically linearly independent. Since $H_u$
is one to one on ${\rm Im}(H_u)$, this is true as well for the
vectors $H_u^{2k+1}(1), k=0,\cdots ,N$, in other words
$${\rm det}(J_{2(m+n+1)})_{0\leq m,n\leq N}\ne 0$$
generically on ${\mathcal M}(N+1)$. Now apply the Cayley-Hamilton
Theorem to $H_u^2$, as we did for the proof of Lemma \ref{levelset}.
For every
   $p=1,\cdots ,N+1,$, we obtain
\begin{equation}\label{JCHN}
J_{2(N+1+p)}=\sum _{j=1}^{N+1}(-1)^{j-1}\sigma _jJ_{2(N+1-j+p)}\ .
\end{equation}
 Solving this linear system, we infer that, locally at generic
 points,
$$\sigma _j=F_j(J_{2k},k=1,\cdots, 2N+2)$$
where $F_j$ is real analytic. Applying again (\ref{JCHN}) for $p=0$,
and observing that $J_0=1-S$ and $\sigma _N\ne 0$ since $H_u$ is one
to one on ${\rm Im}(H_u)$, we obtain, locally at generic points,
$$S=G(J_{2k},k=1,\cdots, 2N+2)$$
where $G$ is real analytic. This implies that the functions $J_{2k},
k=1,\cdots, 2N+2,$ are generically independent on ${\mathcal
M}(N+1)$, which is property (2) for $N+1$. \s The proof that
property (1) implies property (2) is quite similar, so we just
sketch it. First we enlarge $\tilde {\mathcal M}(N)$ as a connected
holomorphic manifold of the same dimension, which contains a dense
open subset of ${\mathcal M}(N)$ as a hypersurface. This can be
realized by considering the manifold $\tilde {\mathcal M}'(N)=\tilde
{\mathcal M}(N)\cup {\mathcal M}(N)\setminus \tilde {\mathcal
M}(N-1)$ which consists of rational functions $u$ of the form
$$u(z)=\frac {A(z)}{B(z)}\ ,$$
with $A\in \C_{N} [z] $, $B\in \C_N [z],$ $B(0)=1,$ $d(A)=N$ or
$d(B)= N$, $A$ and $B$ have no common factors,  and $B(z)\ne 0$ if
$|z|\leq 1$. The coefficient $a$ of $z^N$ in the numerator $A$
defines a holomorphic coordinate on $\tilde {\mathcal M}'(N)$, and
${\mathcal M}(N)$ is defined by the equation $a=0$. Moreover,
$\tilde S=|a|^2$ is a conservation law. Starting from a generic
point $u_0\in {\mathcal M}(N)$, we prove similarly that the
functions $\tilde S, J_{2k}, k=1,\cdots 2N,$ are generically
independent on $\tilde {\mathcal M}'(N)$. Then we infer the generic
independence of $J_{2k}, k=1,\cdots 2N+1,$ by using again the
Cayley-Hamilton theorem for $H_u^2$. \s It is now easy to conclude,
generically on the data in ${\mathcal M}(N)$ or $\tilde {\mathcal
M}(N)$, that the solution of equation (\ref{S}) is quasiperiodic.
Let us sketch the argument for ${\mathcal M}(N)$, for instance. By
Lemma \ref{levelset}, for generic $u_0$ in ${\mathcal M}(N)$, the
level set
$$L_N(u_0):=\{ u\in {\mathcal M}(N): J_{2n}(u)=J_{2n}(u_0), 1\le n\le 2N\} $$
is compact. Moreover, from the generic independence of the functions
$J_{2n}$ combined with the Sard theorem, for generic $u_0\in
{\mathcal M}(N)$, the vector $(J_{2n}(u_0))_{1\le n\le 2N}$ is a
regular value of the mapping
$$u\mapsto (J_{2n}(u))_{1\le n\le 2N}\ .$$
We conclude from standard arguments --- see for instance \cite{Ar},
that, generically on $u_0\in {\mathcal M}(N)$, the level set
$L(u_0)$ is a finite union of $2N$ dimensional Lagrangian tori, on
which the evolution defined by (\ref{S}) is quasiperiodic.
\end{proof}
\section{Traveling waves}\label{TravelingWaves}
We start with some basic definitions. General definitions can be
found in \cite{GSS2}, for example.
 \begin{definition}\label{deftrav}
 A solution $u$ of (\ref{S})  is said to be a traveling wave if there exists $\omega,c\in \R$ such  that $$u(t,z)=e^{-i\omega t}u(0, e^{-ict}z)$$
 for every $t\in \R $. We shall call $\omega $ the {\it pulsation} of $u$, and $c$ the velocity of $u$.
\end{definition}
Equivalently, $u$ is a traveling wave with pulsation $\omega $ and
angular velocity $c$ if and only if it satisfies at time $t=0$ ---
hence at every time--- the following equation,
\begin{equation}\label{trav}
cDu+\omega u=\Pi (|u|^2u)\ .
\end{equation}
In the sequel, a solution $u\in H^{1/2}_+$ of equation (\ref{trav})
will be called as well a traveling wave of pulsation $\omega $ and
of velocity $c$. Notice that equation (\ref{trav}) is variational :
traveling waves of pulsation $\omega $ and velocity $c$ are the
critical points of the functional
$$u\in H^{1/2}_+\mapsto cM(u)+\omega Q(u)-\frac 12 E(u)\ .$$
For example, from Proposition \ref{CaractM(1)}, we know that
elements of ${\mathcal M}(1)$ are characterized as minimizers of
$$u\in H^{1/2}_+\mapsto Q(u)^2+2M(u)Q(u)-E(u)\ ,$$
so that we recover that they are traveling waves with
$$\omega = Q(u)+M(u)\ ,\ c=Q(u)\ .$$

\begin{subsection}{Characterization  of stationary waves}
Stationary waves are  traveling waves with velocity $c$ equal to
$0$. They are particularly easy to characterize.
\begin{proposition}\label{stat}
Let $u_0\in H^{\frac 12}_+\setminus \{ 0\} $. Then $u(t)=e^{-i\omega
t}u_0$ solves $(S)$ if and only if
$$|u_0|^2=\omega \ {\rm a.e.\quad  on}\ \S^1\  ,$$
or equivalently
$$u_0(z)=\alpha \prod _{j=1}^N\frac{z-\overline p_j}{1-p_jz}$$
for some $p_1,\dots ,p_N$ in the unit disc, and $\alpha $ is a
complex number such that $|\alpha |^2=\omega $.
\end{proposition}
\begin{proof} Indeed, $\Pi (|u_0|^2u_0)=\omega u_0$ means
$$|u_0|^2u_0-\omega u_0\, \perp \ L^2_+ $$
which implies $|u_0|^4-\omega |u_0|^2=0$, or $|u_0|^2=\omega $. In
other words, $\varphi :=\omega ^{-1/2}u_0$ is an inner function in
the sense of Beurling. Since $\varphi \in H^{1/2}_+$, the finiteness
of $$(D\varphi |\varphi )=\int _{\S^1}\frac{\varphi '(z)}{\varphi
(z)}\, \frac{dz}{2i\pi} $$ implies, by Rouch\'e's theorem, that
$\varphi $ has only a finite number of zeroes in the unit disc,
therefore is a finite Blaschke product, as claimed.
\end{proof}
As it is well known (see {\it e.g.} \cite{R}, Chapter 17), any inner
function may be written as a product of a Blaschke product and of
$$\exp\left (-\int_0^{2\pi} \frac{e^{it}+z}{e^{it}-z}
d\mu(t)\right)$$ where $\mu$ is a singular measure with respect to
the Lebesgue measure. The simplest cases are
\begin{eqnarray*}u_0(z)=\prod _{j=1}^\infty \frac{z-\overline  {p_j}}{1-p_jz}\ ,\ |p_j|<1,\; \sum _{j=1}^\infty (1-|p_j|)<\infty ,\\
u_0(z)={\rm exp}\left (-\frac{1+z}{1-z}\right )\hskip 1.5cm\ .
\end{eqnarray*}
Let us emphasize that  these particular solutions do not belong to
$H^{1/2}_+$. Hence, these examples show that there exists a larger
family of non smooth solutions of $(S)$, which does not fit with the
existence result of Theorem \ref{Cauchy} and therefore calls for the
construction of a flow map
 on a wider phase space. In view of the $BMO$ conservation law derived from the Lax pair and Nehari's Theorem, a natural candidate for this phase space
 is $BMO_+$. This is a very interesting open question.
 \end{subsection}
\subsection{ Characterization of traveling waves.}
We now focus on the case of a non zero velocity. The main result of
this section is the following.
\begin{theorem}\label{TravWaves}
A function $u\in H^{1/2}_+$ is a traveling wave
 with a velocity $c\in \R^*$  and with a pulsation $\omega\in\R$ if and only if
there exist  non negative integers $N$, $\ell\in\{0,1,\dots,N-1\}$,
and complex numbers $p\in \C$ with $0<|p|<1$  and $\alpha\in \C$,
such that
$$u(z)=\frac{\alpha z^\ell}{1-p^Nz^N}.$$
\end{theorem}
\begin{proof}
We first reformulate the soliton equation (\ref{trav}) in terms of
the Hankel operator $H_u$. Introducing the operator
$$A=D-\frac{1}{c}T_{|u|^2}$$
we observe from (\ref{Rio}) that (\ref{trav}) is equivalent to
\begin{equation}\label{eqop}
AH_u+H_uA+\frac{\omega }{c}H_u+\frac{1}{c}H_u^3=0\ .
\end{equation}
The operator $\tilde A=A+\frac{1}{2c}H_u^2$ is selfadjoint on
$L^2_+$, bounded from below and with a compact resolvent. Therefore
it admits an orthonormal basis of eigenfunctions associated to a
sequence of real eigenvalues tending to $+\infty $. Since
(\ref{eqop}) is equivalent to $$\tilde AH_u+H_u\tilde
A=-\frac{\omega }{c}H_u\ ,$$
 we observe that
$$\tilde A\varphi =\lambda \varphi $$
 yields
$$\tilde  A H_u\varphi =-(\frac{\omega }{c} +\lambda )H_u\varphi $$
and the boundedness of $\tilde A$ from below implies $H_u\varphi =0$
for $\lambda $ large enough. Consequently, $H_u $ has finite rank,
and therefore $u$ is a rational function by the Kronecker theorem
(see appendix 3 for an elementary proof). The main step  is now to
prove the following result.
\begin{proposition}\label{H2U}
There exists $\lambda\in\R$ so that $H^2_u(u)=\lambda u$.
\end{proposition}

Assume this proposition is proved, and let us show how to complete
the proof of Theorem \ref{TravWaves}. We may assume that $1\not\in
{\rm Im}H_u$, otherwise Proposition \ref{H2U} would lead to
$H^2_u(1)=\lambda$ which implies that $|u|^2=\lambda $ and hence
that $u$ is a stationary wave. Denote by $N$ the rank of $H_u$.
Notice that (\ref{eqop}) implies
\begin{equation}\label{comm}
[A,H_u^2]=0
\end{equation}
therefore the range of $H_u^2$ --- which is also the range of
$H_u$--- is invariant through the action of $A$.

As $1\not \in {\rm Im}H_u$, Proposition \ref{H2U} reads
$H^2_u(1)=\lambda P_u(1)$ with $\lambda=\frac Q{1-S}$. Setting
$v=1-P_u(1)$ as in subsection \ref{blaschke}, we have
$$|u|^2=H^2_u(1)+\overline{H^2_u(1)}-Q= \frac{Q}{1-S}(2-v-\overline
v)-Q\ .$$ On the other hand, as $v$ belongs to the kernel of $H_u$,
we have from (\ref{eqop}),
$$H_uA(v)=0\ .$$
But $$A(v)=-\frac 1c H_u^2(1)-AP_u(1)\, \in {\rm Im}(H_u).\ $$ We
conclude that $A(v)=0$, which reads, since $\overline uv$ is
holomorphic from Lemma \ref{ubarh},
$$Dv=\frac{1}{c}|u|^2v\ .$$
>From Rouch\'e 's theorem, we infer
$$Q=Nc.$$
Eventually, we get $$Dv=\frac N{1-S}(2-v-\overline v)v-Nv=\frac{
N(1+S)}{1-S}v- \frac N{1-S}v^2-\frac{NS}{1-S}$$ since $|v|^2=S$.
Notice that the constant $S$ is a particular solution of this
Riccati equation. Solving this equation, we get, for some constant
$B$,
$$v=1-\frac{(1-S)B}{B+z^N}\ .$$
>From this formula, we have
$$H^2_u(1)=\frac{Q}{1-S}(1-v)=\frac {Q}{1-p^Nz^N}$$
for some constant $p$   which is necessarily of modulus  less than
$1$ since $H^2_u(1)$ is holomorphic in the unit disc. It remains to
use that $u$ is solution to the equation
$$cDu+\omega u=H^3_u(1)+uH^2_u(1)-Qu$$ to get that
$$cDu+\omega u=u\left(\frac{QS}{1-S}+\frac {Q}{1-p^Nz^N}\right) .$$
This is an ordinary first order differential equation, which can be
rewritten as
$$D\log(u)=D\log(1-p^Nz^N)^{-1}+N\left (\frac 1{1-S}-\frac\omega Q\right )D\log z\ .$$
By Rouch\'e's theorem, since $u$ is a rational function with no
poles in the unit disc and at most $N-1$ zeroes, we have $$ N\left
(\frac 1{1-S} -\frac \omega Q\right )=\ell \in \{ 0,1,\cdots ,N-1\}
\ .$$ Coming back to the equation on $u$, this proves the claim.
\subsubsection{Proof of Proposition \ref{H2U}}
We now turn to the main step of the proof. Because of (\ref{comm}),
there exists an orthonormal basis of ${\rm Im}(H_u)$ which consists
of common eigenvectors to $A$ and to $H_u^2$. Our strategy is to
describe  precisely the corresponding joint spectrum. Let us
introduce some notation. For  $\gamma>0$, set
$$E_{\lambda,\gamma}=\ker (A-\lambda)\cap \ker(H_u^2-\gamma),$$
and define
$$\Sigma =\{ (\lambda ,\gamma )\in \R\times \R _+^*\, :\,
E_{\lambda,\gamma}\neq \{ 0\} \ \} \ .$$ The following two lemmas
give important information about $\Sigma $. The first one takes
advantage of the relationship with the shift.
\begin{Lemma}\label{shift}
\begin{enumerate}
\item Assume $A\varphi =\lambda \varphi .\\$ If $(\varphi |1)=0$, then
$\varphi =z\psi $ with $A\psi =(\lambda -1)\psi $.\\ If $(z\varphi
|H_u^2(1))=0$, then $A(z\varphi )=(\lambda +1)z\varphi $.
\item Assume $H_u^2\varphi =\gamma \varphi $.\\ If $(\varphi |1)=0$
and $(\varphi |zu)=0$, then $\varphi =z\psi $ with $H_u^2(\psi
)=\gamma \psi $.\\ If $(z\varphi |H_u^2(1))=0$ and $(\varphi |u)=0$,
then $H_u^2(z\varphi )=\gamma z\varphi .$
\end{enumerate}
\end{Lemma}
Lemma \ref{shift} is a straightforward consequence of the following
basic identities :
\begin{eqnarray}\label{identities}
\begin{cases}\begin{aligned}
A(zh)&=zA(h)+zh-\frac 1c(zh|H_u^2(1))\ ,\\ (A(h)|1)&=-\frac 1c(h|H_u^2(1))\ ,\\
H_u^2(zh)&=zH_u^2(h)+(zh|H_u^2(1))-(h|u)zu\ .
\end{aligned}
\end{cases}
\end{eqnarray}
The second lemma specifies the action of $H_u$ on eigenfunctions of
$A$.
\begin{Lemma}\label{Hu} Assume $A\varphi =\lambda
\varphi $ and $H_u^2\varphi =\gamma \varphi $. Then
$$AH_u\varphi =-\left (\lambda +\frac{\omega +\gamma}{c}\right )H_u\varphi
.$$ If, moreover, $(\varphi |1)\ne 0$, then $\gamma =-c\lambda $ and
$$AH_u\varphi =-\frac{\omega}{c}H_u\varphi \ .$$
\end{Lemma}
The first part of Lemma \ref{Hu} is a simple consequence of equation
(\ref{eqop}). The second part follows from the second identity in
(\ref{identities}), which yields
$$\lambda =-\frac{(\varphi |H_u^2(1))}{c(\varphi |1)}=-\frac{(H_u^2\varphi
|1)}{c(\varphi |1)}=-\frac{\gamma}{c}\ .$$ Now we gather the
important facts deduced from the above two lemmas.
\begin{Lemma}\label{lemE}
The following properties hold.
\begin{enumerate}
\item $H_u(E_{\lambda,\gamma})=E_{-(\lambda+\frac{\omega+\gamma}c),\gamma}$.
\item If $E_{\lambda,\gamma}\not \subset 1^\perp$ then $\gamma=-c\lambda$.
\item If $\lambda \neq 1-\frac\omega c$ and $\gamma\neq -c\lambda\ $,
then $E_{\lambda,\gamma}\subset zE_{\lambda-1,\gamma}\ $.
\item If $\lambda \neq 1-\frac\omega c$ and ${\rm dim}(E_{\lambda ,\gamma})\geq 2$, then
$(\lambda -1,\gamma )\in \Sigma \ $.
\end{enumerate}
\end{Lemma}
{\sl Proof.} Lemma \ref{Hu} gives that
$H_u(E_{\lambda,\gamma})\subset
E_{-(\lambda+\frac{\omega+\gamma}c),\gamma}$. For the converse
inclusion, we use the fact that, since $\ker (H_u^2-\gamma)\subset
{\rm Im}H_u^2$ for $\gamma>0$,  any $\varphi\in
E_{-(\lambda+\frac{\omega+\gamma}c),\gamma}$ may be written as
$\varphi=H_u(\psi)$ with $\psi \in {\rm Im}(H_u)$. Since
$H_u^2(\varphi)=\gamma \varphi$ and $H_u$ is one to one on ${\rm
Im}H_u$, we get $H_u(\varphi)=\gamma\psi$ so that $A(\psi)=\frac
1\gamma A(H_u(\varphi))$. We then use Equation (\ref{eqop}) to get
$A(\psi)=\lambda\psi$.

The second assertion is a direct consequence of Lemma \ref{Hu}.

Let us prove the third assertion. Assume $\gamma\neq -c\lambda \  $.
Given $\varphi\in E_{\lambda,\gamma}$, assertion 2 gives
$(\varphi|1)=0$, and Lemma \ref{shift} yields $\varphi=z\psi$ with
$A\psi=(\lambda-1)\psi$. On the other hand,
$(\varphi|zu)=(z\psi|zu)=(\psi|u)=0$ since $u\in \ker (A-\frac\omega
c)$ and $\lambda -1\neq \frac{-\omega }{c}\ $. Hence, by Lemma
\ref{shift}, we have $H^2_u(\psi)=\gamma\psi$ as it is expected.

The proof of the fourth assertion  is a modification of the latter,
based on the following observation : if ${\rm dim}(E_{\lambda
,\gamma})\geq 2$, then $E_{\lambda ,\gamma}\cap 1^\perp \neq \{ 0\}
$. The rest of the proof is unchanged.

\end{proof}

\noindent A consequence is the following description of the joint
spectrum.
\begin{Lemma}\label{PossibleValues}
Given $\gamma >0$, define $\Sigma _\gamma =\{ \lambda \in \R \ :\
(\lambda ,\gamma )\in \Sigma \} \ .$ If $\Sigma _\gamma $ is not
empty, then there exists a nonnegative integer $\ell $ such that one
of the following situations occurs:
\begin{enumerate}
\item Either $\gamma= \omega-(\ell+2)c$ and
$$\Sigma _\gamma =\left \{ 1-\frac \omega c +j\ ,\ j=0,\cdots ,\ell \right \} $$
with the following equalities,
$$E_{-\frac{\omega}c+\ell+1,\omega-(\ell+2)c}= zE_{-\frac{\omega}c+\ell,\omega-(\ell+2)c}= \dots= z^\ell E_{1-\frac{\omega}c,\omega-(\ell+2)c}\ .$$
\item Or $\gamma=\omega+\ell c$ and
$$\Sigma _\gamma =\left \{ -\frac \omega c -j\ ,\ j=0,\cdots ,\ell \right \} $$
with the following equalities,
$$E_{-\frac{\omega}c,\omega+\ell c}=zE_{-\frac{\omega}c-1,\omega+\ell c}= \dots= z^\ell E_{-\frac{\omega}c-\ell,\omega+\ell c}$$
each of the spaces being of dimension $1$.
\end{enumerate}
\end{Lemma}
\begin{proof}
By the third assertion of Lemma \ref{lemE}, if $(\lambda ,\gamma
)\in \Sigma $, then
\begin{enumerate}
\item either $\lambda+\frac{\omega}c$ is an integer $\geq 1$,
\item or $\lambda+\frac\gamma c$ is an integer $\geq 0\ $.
\end{enumerate}
Indeed, otherwise there would exist an infinite sequence of non
trivial eigenspaces
$$E_{\lambda,\gamma}\subset zE_{\lambda-1,\gamma}\subset\dots \subset z^j E_{\lambda-j,\gamma}\subset \dots$$
since for any $j\neq 0$, $\lambda-j\neq 1-\frac{\omega}c$ and
$\gamma\neq-c(\lambda-j)$. This would contradict the boundedness of
$A$ from below. \s
 Applying assertion 1 of Lemma \ref{lemE}, these constraints also
apply to the pair $(\lambda ', \gamma )$ with $$\lambda '=-\lambda
-\frac{\gamma +\omega}{c}\ .$$ This implies
\begin{enumerate}
\item or $\lambda +\frac \gamma c$ is an integer $\leq -1\ $,
\item either $\lambda+\frac{\omega}c$ is an integer $\leq 0$
.
\end{enumerate}
In other words, there exists some nonnegative integer $\ell $ such
that
\begin{enumerate}
\item either $\gamma= \omega-(\ell+2)c$ and
$$\Sigma _\gamma \subset \left \{ 1-\frac \omega c +j\ ,\ j=0,\cdots ,\ell \right \} ,$$
\item or $\gamma=\omega+\ell c$ and
$$\Sigma _\gamma \subset \left \{ -\frac \omega c -j\ ,\ j=0,\cdots ,\ell \right \} \ .$$
\end{enumerate}
Assume now that, say $\gamma= \omega-(\ell+2)c$. Applying assertion
3 of Lemma \ref{lemE}, we obtain, for some $k\in \{ 0,\cdots,
\ell\}$,
$$\{ 0\} \neq E_{-\frac{\omega}c+k+1,\omega-(\ell+2)c}\subset zE_{-\frac{\omega}c+k,\omega-(\ell+2)c}\subset
 \dots\subset z^{k}E_{1-\frac{\omega}c,\omega-(\ell+2)c}\ .$$
Applying assertion 1 of Lemma \ref{lemE}, and again assertion 3, we
also have
$$H_u(E_{1-\frac{\omega}c,\omega-(\ell+2)c})=E_{-\frac{\omega}c+\ell+1,\omega-(\ell+2)c}\subset
  \dots\subset z^{\ell -k} E_{-\frac{\omega}c+k+1,\omega-(\ell+2)c}\ .$$
Consequently, we have the claimed equality by a dimension argument.
\s The same procedure applies to the case $\gamma =\omega +\ell c$.
Moreover, by assertion 4 of Lemma \ref{lemE}, we know that the
dimension of $E_{-\frac{\omega}c-\ell,\omega+\ell c}$ is at most
$1$, hence exactly $1$, which completes the proof.
\end{proof}
\begin{proof}
We now turn to the proof of Proposition \ref{H2U} itself. We argue
by contradiction and assume that $H^2_u(u)$ and $u$ are independent
so that the eigenvalue $-\frac\omega c$ of $A$ is not simple. As a
first consequence of the fourth assertion of Lemma \ref{lemE}, the
minimal eigenvalue of $A$ on ${\rm Im}(H_u)$ is necessarily simple.
By Lemma \ref{PossibleValues}, since $-\frac\omega c$ is an
eigenvalue of multiplicity at least $2$, this minimal eigenvalue is
necessarily of the form $\lambda_{min}=-\frac{\omega}c-j$ for some
positive integer $j$. Again, by Lemma \ref{PossibleValues}, we
therefore have $\ker(A+\frac\omega c)\cap{\rm Im}H_u=\oplus_{k\in
K}E_{-\frac{\omega}c,\omega+kc}$ where $K$ is a finite subset of
$\{0,\dots,j\}$ containing at least $j$ and another integer.
Furthermore, all the spaces $E_{-\frac{\omega}c,\omega+kc}\ $, $k\in
K$, have dimension $1$.  \s We are going to prove that $K$ has
exactly two elements. Our strategy is based on the following
observation, which is a direct consequence of Lemma \ref{shift} : if
$\varphi \in \ker(A+\frac\omega c)$ satisfies  $(z\varphi
|H_u^2(1))=0$, then $z\varphi$ belongs to $\ker(A-1+\frac\omega c)$.
Consequently,
$$|K|={\rm dim}\left (\ker (A+\frac \omega c)\cap {\rm Im}(H_u)\right )\leq
1+{\rm dim}(\mathcal N)\ ,$$ where
$$\mathcal N:=\ker (A+\frac \omega c -1)\cap z\, \left (\ker (A+\frac \omega c)\cap
{\rm Im}(H_u)\right )$$ and, if we prove that $\mathcal N$ is at
most one dimensional, we will conclude that $|K|=2$. \s
 As a first step, we are going to study the auxiliary space
$\ker(A-1+\frac\omega c)\cap {\rm Im}(H_u)$. By Lemma
\ref{PossibleValues}, this space is the direct sum of spaces $
E_{1-\frac{\omega}c,\gamma }$, where $\gamma $ describes a set of
positive values included in $\{ \omega -(\ell +2)c, \ell =0,1,\cdots
\} $. In view of assertion 2 of Lemma \ref{lemE}, elements $\psi $
of
 $ E_{1-\frac{\omega}c,\gamma }$ satisfy
$(\psi|1)=0$, hence we can write $\psi=z\varphi$ with $\varphi\in
\ker(A+\frac\omega c)$, because of assertion 1 of Lemma \ref{shift}.
Moreover, using the third formula of (\ref{identities}), the
equation $H_u^2\psi =\gamma \psi $ reads $(z\varphi |H_u(u))=0$ and
\begin{equation}\label{ConditionH2phi}
H^2_u\varphi=\gamma \varphi +( \varphi|u)\, u\ ,\ .
\end{equation}
hence $\varphi \in {\rm Im}(H_u)\cap \ker (A+\frac \omega c)\ $. Let
us compute the characteristic polynomial of the eigenvalue problem
(\ref{ConditionH2phi}) on ${\rm Im}(H_u)\cap \ker (A+\frac \omega
c)$. Let $\{\varphi_k\}_{k\in K}$ be an orthonormal basis of
$\ker(A+\frac\omega c)\cap{\rm Im}H_u=\oplus_{k\in K}E_{-\frac
\omega c,\omega+kc}$, with $\varphi_k\in
E_{-\frac{\omega}c,\omega+kc}$ for any $k\in K$. We write
$$\varphi=\sum_{k\in K}\alpha_k\varphi_k,\; u=\sum_{k\in K}\beta_k\varphi_k$$
Computing both sides of (\ref{ConditionH2phi}) in coordinates, we
get
$$\sum_{k\in K}\alpha_k(\omega+ck)\varphi_k=\sum_{k\in K}(\gamma\alpha_k+\beta_k(\sum _{k'\in K}\alpha_{k'}\overline{\beta_{k'}}))\varphi_k$$
so that the $\alpha_k$'s have to satisfy the following system
$$\alpha_k(\omega +ck -\gamma)=\beta_k\sum _{k'\in K}\alpha_{k'}\overline{\beta_{k'}}\ .$$
 The
characteristic polynomial is the determinant of this system, namely
\begin{equation}\label{charac}
P(\gamma )=\prod_{k\in K}(\omega +kc-\gamma )\left(1-\sum_{k\in
K}\frac{|\beta_k|^2}{\omega +kc-\gamma}\right)\ .
\end{equation}
Plugging the additional information $\gamma =\omega -(\ell +2)c$ for
some nonnegative integer $\ell $, the equation is then equivalent to
\begin{equation}\label{EqDet}\sum_{k\in K}\frac{|\beta_k|^2}{(k+\ell +2)c}=1\
,
\end{equation}
which admits a unique simple  solution in $\ell $ if $c>0$, and no
solution if $c<0$. Hence $\ker(A-1+\frac\omega c)\cap {\rm Im}(H_u)$
is $\{ 0\}$ if $c<0$, and is at most one dimensional if $c>0$. \s
Next we distinguish two cases. \s {\sl First case: $1\notin {\rm Im}
H_u$.} Then  the kernel of $H_u$ is $bL^2_+$, where $b$ is a finite
Blaschke product  with $b(0)\ne 0$. We infer that $$z{\rm
Im}(H_u)\cap \ker H_u=\{ 0\} .$$ Indeed, if $zH_u(\varphi )=bh$,
then $h$ is divisible by $z$ and thus $H_u(\varphi )\in bL^2_+=\ker
H_u$, hence $H_u(\varphi )=0$. Now we consider the orthogonal
projection onto ${\rm Im}H_u$ restricted to $\mathcal N$. The kernel
of this linear mapping is contained into $z{\rm Im}(H_u)\cap \ker
H_u$, therefore this mapping is one to one. Since its image is
contained into $\ker(A-1+\frac\omega c)\cap {\rm Im}(H_u)$, which is
at most one dimensional, $\mathcal N$ is at most one dimensional. We
conclude that $|K|=2$. We notice that, in this case, we have proved
that $\ker(A-1+\frac\omega c)\cap {\rm Im}(H_u)$ is exactly one
-dimensional.
 \s {\sl Second case: $1\in {\rm
Im} H_u$.} In this case, we shall determine $\ker(A-1+\frac\omega
c)$ itself. Let us  make some preliminary remarks. Recall from
Proposition \ref{w} that the solution $w\in {\rm Im}H_u$ of
$$H_u(w)=1$$
satisfies $zw=Cb$ for some constant $C$ where $b$ is a Blaschke
product of degree $N$, the Beurling generator of $\ker H_u$. From
(\ref{eqop}),
$$A(w)+\frac{\omega}{c}w=0$$
or, since $\overline uw$ is holomorphic by Lemma \ref{ubarw},
$$Dw+\frac{\omega}{c}w=\frac 1c|u|^2w\ ,\ Db+\left (\frac{\omega}{c}-1\right )b=\frac 1c|u|^2b$$
whence, again by Rouch\'e 's theorem,
\begin{equation}\label{rouche}
Q=(N-1)c+\omega .
\end{equation}
Observe that the above equation on $b$ means that $$b\in \ker
(A-1+\frac \omega c)\ .$$ Moreover, $\ker (A-1+\frac \omega c)\cap
\ker H_u$ consists of functions $bh$ satisfying
$$D(bh)-\frac 1c|u|^2bh+\left (\frac \omega c-1\right )bh=0\ , $$
or $Dh=0$. Hence
$$\ker (A-1+\frac \omega c)\cap \ker H_u=\C b\ .$$
It remains to describe $\ker(A-1+\frac\omega c)\cap {\rm Im}(H_u)$.
We already know that this space is $\{ 0\} $ if $c<0$.  To study the
case $c>0$, we return to equation (\ref{ConditionH2phi}). We observe
that $\varphi =w$ is a solution of this equation with $\gamma =0$,
since $H^2_u(w)=u=(w|u) u$. Moreover, the characteristic polynomial
$P(\gamma )$ given by (\ref{charac}) admits a unique zero in the
interval $(-\infty ,\min _{k\in K}(\omega +kc))$. Since this
interval contains all the values $\omega -(\ell +2)c\ ,\ \ell
=0,1,\cdots \ $, and $0$ --- indeed $\omega +kc , k\in K$, is an
eigenvalue of $H_u^2$ on ${\rm Im}(H_u)$, hence is positive ---
 we
conclude that $$\ker(A-1+\frac\omega c)\cap {\rm Im}(H_u)=\{ 0\}\
.$$ Therefore $\ker(A-1+\frac\omega c)=\C b$, so that $\mathcal N$
is at most one dimensional and $|K|=2$. \s We can finally write
$$\ker(A+\frac\omega c)\cap {\rm Im}H_u=E_{-\frac \omega
c,\omega+jc}\oplus E_{-\frac \omega c,\omega+kc}$$ with $0\leq k<j$.
\s As a final step, we are going to get a contradiction implied by
this two-dimensionality. \s
 We first consider the
case when $1\in{\rm Im}(H_u)$. Let us apply the Cayley-Hamilton
theorem to $H_u^2$ on the two-dimensional space $\ker(A+\frac \omega
c)$. We obtain
$$H_u^4(u)- (2\omega+(j+k)c)H_u^2(u)+(\omega+jc)(\omega+kc)u=0\ .$$
Since $1\in {\rm Im}(H_u)$, this implies
$$H_u^4(1)-
(2\omega+(j+k)c)H_u^2(1)+(\omega+jc)(\omega+kc)=0\ ,$$ and, taking
the scalar product with $1$,
$$J_4-(2\omega+(j+k)c)Q+(\omega+jc)(\omega+kc)=0\ .$$
Using that $Q=(N-1)c+\omega$ by (\ref{rouche}), we get
\begin{eqnarray*}
J_4-Q^2&=& (2\omega+(j+k)c-(N-1)c-\omega)((N-1)c+\omega)-(\omega+jc)(\omega+kc)\\
&=&-c^2(N-1+jk)<0 \ .
\end{eqnarray*}
This fact is in contradiction with the Cauchy-Schwarz inequality,
$$Q^2=|(H_u^2(1)|1)|^2\leq \| H_u^2(1)\| ^2=(H_u^4(1)|1)=J_4\ .$$
It remains to consider the case $1\notin {\rm Im}(H_u)$. Again, we
are going to contradict the Cauchy-Schwarz inequality. First, we use
the Cayley-Hamilton theorem as before,
$$H_u^4(1)-
(2\omega+(j+k)c)H_u^2(1)+(\omega+jc)(\omega+kc)P_u(1)=0\ ,$$ which
yields to
$$J_4-(2\omega+(j+k)c)Q+(\omega+jc)(\omega+kc)(1-S)=0\ $$
and
\begin{equation}\label{degre2}J_4(1-S)-Q^2=-\, \frac{(J_4-Q(\omega
+jc))(J_4-Q(\omega +kc))}{(\omega+jc)(\omega+kc)}\ . \end{equation}
The Cauchy--Schwarz inequality
$$Q^2=|(H_u^2(1)|P_u(1))|^2\leq \| H_u^2(1)\| ^2\| P_u(1)\| ^2=J_4(1-S)$$
implies that the left hand side of (\ref{degre2}) is nonnegative.
Therefore, remembering that $\omega +jc$ and $\omega +kc$ are
positive as eigenvalues of $H_u^2$ on ${\rm Im}(H_u)$, we shall
obtain a contradiction if we show that
\begin{equation}\label{contrad}
J_4\, >\, Q(\omega +jc)\ .
\end{equation}
Let us prove (\ref{contrad}). Recall that $c>0$, since $Q=Nc$. Apply
Lemma \ref{PossibleValues}. If $\gamma
>0$ is an eigenvalue of $H_u^2$, either $\gamma =\omega +\ell c$
with $\ell \ge 0$, and $E_{-\frac \omega c, \gamma }\ne \{ 0\} $,
and this implies $\ell \in \{ j,k\} $ ; or $\gamma =\omega -(\ell
+2)c$, $\ell \ge 0$, and $E_{1-\frac \omega c, \gamma }\ne \{ 0\}$.
In this case, we have already seen that $\ker(A-1+\frac\omega c)\cap
{\rm Im}(H_u)$ is one dimensional, which means that $\ell $ is
uniquely determined and $E_{1-\frac \omega c, \gamma }$ is
one-dimensional. We infer the following decomposition, where all the
spaces $E_{\lambda ,\gamma}$ are one--dimensional,
\begin{eqnarray*}
 {\rm Im}(H_u)&=&E_1\oplus E_2\oplus E_3\ ,\\
 E_1&=&\oplus _{j'=0}^{j}E_{-\frac{\omega }{c}-j',\omega
 +jc}\, ,\\ E_2&=&\oplus _{k'=0}^k E_{-\frac{\omega }{c}-k',\omega
 +kc}\, ,\\ E_3&=&\oplus _{\ell '=0}^\ell E_{1-\frac{\omega }{c}+\ell ',\omega
 -(\ell +2)c}\ .
 \end{eqnarray*}
 Consequently,  $N=j+k+\ell+3$ and
\begin{eqnarray*}{\rm Tr}(H^2_u)&=&(j+1)(\omega+jc)+(k+1)(\omega +kc) +(\ell +1)(\omega -(\ell +2)c)\\
&=& N\omega +c[j(j+1)+k(k+1)-(\ell +1)(\ell +2)].
\end{eqnarray*}
On the other hand, ${\rm Tr}(H^2_u)=M+Q=M+Nc$, and, taking the
scalar product of $u$ with both sides of  the soliton equation
(\ref{trav}), we have,
$$M+\frac{\omega }{c}Q=\frac{1}{c}(2J_4-Q^2)\ .$$
Using the identity $Q=Nc$, we infer
$$2J_4=Mc+N\omega c+N^2c^2\ ,$$
and, using the above expression of $M$ provided by the trace of
$H_u^2$,
$$2J_4=2N\omega c+c^2(N^2+j(j+1)+k(k+1)-(\ell +1)(\ell +2)-N)\ .$$
Consequently,
\begin{eqnarray*}
\begin{aligned}2(J_4-Q(\omega +jc))&=2J_4-2N\omega c-2Njc^2\hfill
\\&=c^2(N^2+j(j+1)+k(k+1)-(\ell +1)(\ell +2)-N(2j+1))\\
&=2c^2(k+1)(k+\ell +2)>0 \end{aligned}
\end{eqnarray*} as can be shown by a
straightforward calculation. This proves (\ref{contrad}) and yields
the contradiction, completing the proof of Theorem \ref{TravWaves}.
\end{proof}

\section{Appendices}
\subsection{Appendix 1: The Brezis Gallou\"et estimate}

We recall a simple proof of the estimate
$$\| u\| _{L^\infty }\leq C_s\| u\| _{H^{1/2}}\left [\log \left (2+\frac{\| u\| _{H^s}}{\| u\| _{H^{1/2}}}\right )\right ]^{\frac 12}\ .$$
By Fourier expansion, one has, for any $N\in \N$
\begin{eqnarray*}
&&\|u\|_{L^\infty}\le \sum |\hat u(k)|\\
&=&\sum_{|k|\le N}(1+|k|)^{1/2}\frac{|\hat
u(k)|}{(1+|k|)^{1/2}}+\sum_{|k|\ge N+1}(1+|k|)^{s}\frac{|\hat
u(k)|}{(1+|k|)^{s}}\cr &\le &\|u\|_{H^{1/2}}\times
\left(\sum_{|k|\le N} \frac 1{1+|k|}\right)^{1/2} +\|u\|_{H^s}\times
\left(\sum_{|k|\ge N+1} \frac 1{(1+|k|)^{2s}}\right)^{1/2}\cr &\le
&C\left(\|u\|_{H^{1/2}}\log(N+1)^{1/2}+\|u\|_{H^s}
N^{-s+1/2}\right).
\end{eqnarray*}
The result follows by taking the minimum over $N$.

\subsection{Appendix 2: A Trudinger-type estimate.}

Let us prove the estimate \begin{equation}\label{Ebis} \forall p<
\infty\, , \, \| u\| _{L^p}\leq C\, \sqrt {p}\, \| u\| _{H^{1/2}}\
.\end{equation} It follows from a Marcinkiewicz type argument.
 Assume $\|u\|_{H^{1/2}}=1$. Write, for any $p>2$,
$$\|u\|_{L^p}^p=p\int_0^\infty t^{p-1}\sigma(\{x,\; |u(x)|\ge t\})dt$$
and decompose $u=u_{>\lambda}+u_{<\lambda}$ where
$u_{<\lambda}=\sum_{|k|\le \lambda}\hat u(k)e^{ik\theta}$. Choose
$\lambda=\lambda_t$ so that $\|u_{<\lambda}\|_\infty\le t/2$. More
precisely, since
\begin{eqnarray*}
\|u_{<\lambda}\|_\infty&\le& \sum_{|k|\le \lambda}|\hat u(k)|\\
&\lesssim&  \left(\sum_{|k|\le \lambda }(|k|^2+1)^{1/2}|\hat u(k)|^2\right)^{1/2}\times[\log(\lambda+1 )]^{1/2}\\
&\lesssim&
\|u\|_{H^{1/2}}[\log(\lambda+1)]^{1/2}=c[\log(\lambda+1)]^{1/2} ,
\end{eqnarray*}
we can choose $\lambda$ so that $c[\log(\lambda+1)]^{1/2}=\frac t2$.
With this choice, we get
\begin{eqnarray*}
\|u\|_{L^p}^p&\le& p\int_0^\infty t^{p-1}\sigma(\{x,\; |u_{>\lambda_t}(x)|\ge t/2\})dt\\
&\le& p\int_0^\infty t^{p-3}\|u_{>\lambda_t} \|_2^2dt\le p\int_0^\infty t^{p-3}\sum_{|k|\ge\lambda_t}|\hat u(k)|^2dt\\
&\le &p\sum_{k}\left(\int_0^{2\log(|k|+1)^{1/2}} t^{p-3} dt\right)|\hat u(k)|^2\\
&\le &\frac p{p-2}\sum_k (\log(|k|+1))^{(p-2)/2} |\hat u(k)|^2.
\end{eqnarray*}

Eventually, we use that $(\log(|k|+1))^{\ell}\lesssim \ell!(|
k|+1)\lesssim \ell^\ell(|k|^2+1)^{1/2}$.
 It gives the expected constant proportional to $p^{1/2}$ in (\ref{Ebis}).

\subsection{Appendix 3: An elementary proof of the Kronecker Theorem.}

Let $u\in BMO_+(\S^1)$ so that the Hankel operator $H_u$ is well
defined as a bounded operator on $L^2_+(\S^1)$. Since $H_u$ is $\C$
-antilinear, the range of $H_u$ is a complex vector space.

\begin{proposition}
The function $u$ belongs to ${\mathcal M}(N)$ if and only if the
Hankel operator $H_u$ has (complex) rank $N$. Moreover, if
$$B(z)=\prod _{j=1}^N (1-p_jz)$$
is the denominator of $u$,  the image of $H_u$ is the vector space
generated by
$$\frac{1}{(1-pz)^m}$$
for $0< |p |<1, 1\leq m\leq m_p$, or of the form
$$z^m, 0\leq m\leq m_0-1\ ,$$
where $m_p$ is the number of occurrences of $p$ in the list
$p_1,\cdots ,p_N$.
\end{proposition}
\begin{proof} The proof is based on the following two observations.

i) If $u\in {\mathcal M}(N)$, then ${\rm rk}(H_u)\leq N$.

ii) If ${\rm rk}(H_u)=N$, then $u\in {\mathcal M}(N)$.

Let us first prove i). If $u\in {\mathcal M}(N)$, then one can write
$u$ as a linear combination of functions of the form
$$\frac{1}{(1-pz)^m}$$
for $0< |p |<1, 1\leq m\leq m_p$, or of the form
$$z^m, 0\leq m\leq m_0-1\ ,$$
which we shall associate to $p=0$, with the following degree
condition,
$$\sum _{p}m_p=N\ .$$
Indeed, either the denominator of $u$ is of degree $N$, and this
corresponds to the fact that all the $p$'s are different from $0$,
and the above identity reflects the degree of the denominator ; or
the denominator has degree $<N$, and then the numerator should be of
degree exactly $N-1$ ; therefore the decomposition of $u$ into
elementary fractions involves a polynomial function of degree
$m_0-1\geq 0$, and the above identity reflects the degree of the
numerator $+1$. Now we recall that
$$\widehat {H_u(h)}(k)=\sum _{\ell \geq 0}\hat u(k+\ell )\overline {\hat h(\ell )}\ .$$
In view of the decomposition of $u$, we observe that the sequence
$(\hat u(k))_{k\geq 0}$ is a linear combination of the following
sequences,
$$k^{m-1}p^k, 1\leq m\leq m_p\ ,$$
for $p\ne 0$, and
$$\delta _{km}, 0\leq m\leq m_0-1\ ,$$
for $p=0$. This implies that all the sequences $(\widehat
{H_u(h)}(k))_{k\geq 0}$ have the same property, and therefore the
range of $H_u$ is included into the space $V$  of linear
combinations of
$$\frac{1}{(1-pz)^m}, 1\leq m\leq m_p, 0< |p |<1\ ;\ z^m, 0\leq m\leq m_0-1\ .$$
This implies that ${\rm rk}(H_u)\leq N$.

We now proceed to the proof of property ii). We know that $H_u$ is a
symmetric operator of real rank $2N$. Restricting $H_u$ to its
range, which is a complex vector space of dimension $N$ and is the
orthogonal of ${\rm Ker}(H_u)$ (for both real scalar product and
hermitian scalar product), we can find a real orthonormal basis of
eigenvectors of $H_u$. Moreover, since $H_u$ is antilinear, we
observe that, if $H_u(v)=\lambda v$, then $H_u(iv)=-i\lambda v$.
Therefore we may assume that the above real orthonormal basis of
${\rm Im}(H_u)$ has the special form
$$v_1,iv_1,v_2,iv_2,\dots, v_N,iv_N\ ,$$
and that $H_u(v_j)=\lambda _jv_j$ with some $\lambda _j>0$. Defining
$w_j:=\sqrt {\lambda _j}v_j$, we obtain the following expression for
$H_u$,
$$H_u(h)=\sum _{j=1}^N(w_j|h)_{L^2}\, w_j\ ,$$
or equivalently,
$$\hat u(k+\ell )=\sum _{j=1}^N\hat w_j(k)\hat w_j(\ell )\ ,$$
for all $k\geq 0,\ell \geq 0$. Now the matrix $(\hat w_j(\ell
))_{1\leq j\leq N,0\leq \ell\leq N}$ has rank at most $N$, therefore
there exists $(c_0,c_1,\dots ,c_N)\ne (0,\dots ,0)$ in $\C ^{N+1}$
such that
$$\sum _{\ell =0}^N c_\ell \, \hat w_j(\ell )=0$$
for every $j=1,\dots ,N$. This implies that
$$\sum _{\ell =0}^N c_\ell \, \hat u(k+\ell )=0\ ,$$
for every $k\geq 0$. We then introduce the polynomial
$$P(X)=\sum _{\ell =0}^N c_\ell X^\ell \ .$$
Let
$${\mathcal P}=\{ p \in \C , P(p)=0\} $$
and $m_p\geq 1$ denotes the multiplicity of $p\in {\mathcal P}$.
Then the theory of linear recurrent sequences  implies that the
sequence $(\hat u(k))_{k\geq 0}$ is a linear combination of the
following sequences,
$$k^{m-1}p^k, 1\leq m\leq m_p\ ,$$
for $p\ne 0$, and
$$\delta _{km}, 0\leq m\leq m_0-1\ ,$$
for $p=0$. In other words, $u$ is a linear combination of the
following functions,
$$\frac{1}{(1-pz)^m}\, ,\,  1\leq m\leq m_p\, ,\,  0< |p |<1\ ;\ z^m, 0\leq m\leq m_0-1\ .$$
Since $\sum _pm_p\leq N$, this implies that $u\in {\mathcal M}(N')$
for some $N'\leq N$. However, if $N'<N$, assertion i) would imply
${\rm rk}(H_u)\leq N'$, which contradicts the assumption. Therefore
$N'=N$, and ii) is proved.

Finally, in view of ii), i) is strengthened into\\
i)' If $u\in {\mathcal M}(N)$, then ${\rm rk}(H_u)=N$. \s Moreover,
the inclusion of the range of $H_u$ into the space $V$ becomes an
equality, which is the claim.

This completes the proof.
\end{proof}


\begin{thebibliography}{MTW1}

\bibitem{Ar} Arnold, V.I., {\em Mathematical Methods of Classical
Mechanics}, Springer, New York, 1978.

\bibitem{BGX} Bahouri, H., G\'erard, P., Xu, C.J.: {\em Espaces de Besov et estimations de
Strichartz g\'en\'eralis\'ees sur le groupe de Heisenberg}, Journal
d'Analyse Math\'ematique, 82, 93-118 (2000).

\bibitem{BSV} Birnir, B.,  Kenig, C.,  Ponce, G.,  Svansted, N.,
Vega, L.: {\em On the ill-posedness of the IVP for the generalized
KdV and nonlinear Schr\"odinger equation.} J. London Math. Soc. 53,
551-559 (1996).

\bibitem{B} Bourgain, J.: {\em Refinements of Strichartz' inequality and applications to 2D NLS with critical nonlinearity},
IMRN, 5, 253-283 (1998).

\bibitem{BG} Brezis, H., Gallou\"et, T.: {\em Nonlinear Schr\"odinger
evolution equations.} Nonlinear Anal. 4, 677--681 (1980).

\bibitem{BGT0} Burq, N., G\'erard, P.,  Tzvetkov, N.: {\em
Strichartz inequalities and the nonlinear Schr\"odinger equation on
compact manifolds} Amer. J. Math. 126, 569--605 (2004).

\bibitem{BGT1} Burq, N., G\'erard, P.,  Tzvetkov, N.:  {\em An instability property of the nonlinear {S}chr\"odinger equation on
{$S\sp d$}.}  Math. Res. Lett., 9,  323--335 (2002).

\bibitem{BGT2} Burq, N., G\'erard, P., Tzvetkov, N.:
{\em Bilinear eigenfunction estimates and the nonlinear
Schr\"odinger equation on surfaces.} Invent. math. 159, 187-223
(2005)


\bibitem{BGT3} Burq, N., G\'erard, P., Tzvetkov, N.:
{\em Multilinear eigenfunction estimates and global existence for
the three dimensional nonlinear Schr\"odinger equations.}
 Ann. Scient. \'Ec. Norm. Sup. 38, 255--301 (2005).

\bibitem{BGTstrauss} Burq, N., G\'erard, P., Tzvetkov, N. : {\em High
frequency solutions of the nonlinear Schr\"odinger equation on
surfaces.} Quart. Appl. Math., to appear, 2009.

\bibitem{CKSTT} Colliander, J., Keel, M., Staffilani, G., Takaoka,
H., Tao, T., : {\em Weakly turbulent solutions for the cubic
defocusing nonlinear Schr\"odinger equation}, preprint, 2008, arXiv:
08081742v2 [math.AP].

\bibitem{PG} G\'erard, P.: {\em Nonlinear Schr\"odinger equations in
inhomogeneous media: wellposednes and illposedness results.}
Proceedings of the International Congress of Mathematicians, Madrid,
Spain, 2006, European Mathematical Society.

\bibitem{GG} G\'erard, P., Grellier, S., {\em L'\'equation de Szeg\"o cubique.} S\'eminaire X-EDP, 20 octobre 2008,
\'Ecole Polytechnique, Palaiseau.

\bibitem{GSS2} Grillakis, M., Shatah, J., Strauss, W., {\em Stability theory of solitary waves in the presence of symmetry.
II.}
  J. Funct. Anal.  94 ,  308--348 (1990).

\bibitem{KP} Kappeler, T., P\"oschel, J. : {\em KdV \& KAM}, A
Series of Modern Surveys in Mathematics, vol. 45, Springer-Verlag,
2003.

\bibitem{Kr} Kronecker, L. : {\em Zur Theorie der Elimination einer
Variablen aus zwei algebraische Gleischungen} Montasber. K\"onigl.
Preussischen Akad. Wies. (Berlin), 535-600 (1881). Reprinted in {\em
mathematische Werke}, vol. 2, 113--192, Chelsea, 1968.

\bibitem{K} Kuksin, S. B.: {\em  Analysis of Hamiltonian PDEs.}
 Oxford Lecture Series in Mathematics and its Applications, 19.
 Oxford University Press, Oxford, 2000.

\bibitem{L} Lax, P. : {\em Integrals of Nonlinear equations of
Evolution and Solitary Waves}, Comm. Pure and Applied Math. 21,
467-490 (1968).

\bibitem{L2} Lax, P. : {\em Periodic solutions of the the KdV equation.}  Comm. Pure Appl. Math.  28 , 141--188
(1975).

\bibitem{Ne} Nehari, Z. : {\em On bounded bilinear forms.} Ann.
Math. 65, 153--162 (1957).

\bibitem{N} Nikolskii, N.~K. : {\em Operators, functions, and systems: an easy reading. Vol. 1. Hardy, Hankel, and
Toeplitz.} Translated from the French by Andreas Hartmann.
Mathematical Surveys and Monographs, 92. American Mathematical
Society, Providence, RI, 2002.

\bibitem{Nier} Nier, F. : {\em Bose-Einstein condensates in the Lowest
Landau Levl : Hamiltonian Dynamics.}  Rev. Math. Phys. 19 , 101--130
(2007).

\bibitem{O} Ogawa, T. : {\em A proof of Trudinger's inequality and its application to nonlinear Schr\"odinger
equations.}
  Nonlinear Anal.  14,   765--769 (1990).

\bibitem{P} Peller, V.~V.: {\em Hankel operators and their applications}. Springer Monographs in Mathematics.
 Springer-Verlag, New York, 2003.

\bibitem{R} Rudin, W.: {\em Real and Complex Analysis}, Mac Graw
Hill, Second edition, 1980.

\bibitem{T} Tzvetkov, N.: {\em A la fronti\`ere entre EDP semi et
quasi lin\'eaires}, M\'emoire d'habilitation \`a diriger les
recherches, Universit\'e Paris-Sud, Orsay, 2003.

\bibitem{V} Vladimirov, M. V.: {\em On the solvability of a mixed problem for a nonlinear equation of Schr\"odinger
type.} Sov. math. Dokl. 29, 281-284 (1984).

\bibitem{W1} Weinstein, M. {\em Nonlinear Schr\"odinger equations and sharp interpolation estimates.}
  Comm. Math. Phys.  87   567--576 (1982/83).

\bibitem{W2} Weinstein, M.: {\em Lyapunov stability of ground states of nonlinear dispersive evolution equations.}
  Comm. Pure Appl. Math.  39 , 51--67 (1986).

\bibitem{Y} Yudovich, V. I.: {\em Non-stationary flows of an ideal incompressible fluid.}
 (Russian)  Z. Vycisl. Mat. i Mat. Fiz.  3,  1032--1066 (1963).

\bibitem{Z} Zakharov, V. E., Shabat, A. B.: {\em Exact theory of two-dimensional self-focusing and
one-dimensional self-modulation of waves in nonlinear media.} Soviet
Physics JETP  34  (1972), no. 1, 62--69.
\end{thebibliography}
\end{document}